\documentclass[a4paper,10pt,oneside,reqno]{amsart}

\usepackage[utf8]{inputenc}
\usepackage[hidelinks]{hyperref}

\usepackage[DIV=9]{typearea}

\usepackage{microtype}
\usepackage{amssymb,amsmath,amsopn,amsxtra,amsthm,amsfonts}

\let\originallhook\lhook
\usepackage{MnSymbol} 
\newcommand{\lhook}{\mathrel{\raise.018ex\hbox{$\originallhook$}}}

\usepackage{xr}
\usepackage[mathcal]{euscript}
\let\mathscr\mathcal
\usepackage[bb=px]{mathalfa}

\usepackage{enumitem}
\setlist[enumerate,1]{label={(\arabic*)},itemsep=\parskip} 
\setlist[itemize,1]{itemsep=\parskip} 
\newlist{thmlist}{enumerate}{2}
\setlist[thmlist,1]{label={\em(\roman*)},ref={(\roman*)},%
  itemsep=\parskip,leftmargin=*,align=left}
\setlist[thmlist,2]{label={\em(\alph*)},ref={(\alph*)},%
  itemsep=\parskip,leftmargin=*,align=left,topsep=0.1cm}
\newlist{defnlist}{enumerate}{2}
\setlist[defnlist,1]{label={(\roman*)},ref={(\roman*)},itemsep=\parskip,%
  leftmargin=*,align=left}
\setlist[defnlist,2]{label={(\alph*)},ref={(\alph*)},itemsep=\parskip,%
  leftmargin=*,align=left,topsep=0.1cm}

\newtheoremstyle{plain}
  {.5\baselineskip\@plus.2\baselineskip\@minus.2\baselineskip}
  {.5\baselineskip\@plus.2\baselineskip\@minus.2\baselineskip\@plus.5em}
  {\slshape}
  {}
  {\bfseries}
  {.}
  { }
  {}

\newtheoremstyle{definition}
  {.5\baselineskip\@plus.2\baselineskip\@minus.2\baselineskip}
  {0.8\baselineskip\@plus.2\baselineskip\@minus.2\baselineskip\@plus.5em}
  {}
  {}
  {\bfseries}
  {.}
  { }
  {}

\theoremstyle{plain}

\newtheorem{thm}[subsubsection]{Theorem}
\newtheorem{thmX}{Theorem}

\newtheorem{cor}[subsubsection]{Corollary}
\newtheorem{lem}[subsubsection]{Lemma}
\newtheorem{prop}[subsubsection]{Proposition}

\theoremstyle{definition}
\newtheorem{constr}[subsubsection]{Construction}
\newtheorem{rem}[subsubsection]{Remark}
\newtheorem{defn}[subsubsection]{Definition}

\newtheorem{exam}[subsubsection]{Example}
\newtheorem{notat}[subsubsection]{Notation}

\newcommand{\constrref}[1]{Construction~\ref{#1}}
\newcommand{\thmref}[1]{Theorem~\ref{#1}}

\newcommand{\secref}[1]{Sect.~\ref{#1}}
\newcommand{\ssecref}[1]{Subsect.~\ref{#1}}

\newcommand{\lemref}[1]{Lemma~\ref{#1}}
\newcommand{\propref}[1]{Proposition~\ref{#1}}
\newcommand{\corref}[1]{Corollary~\ref{#1}}

\newcommand{\remref}[1]{Remark~\ref{#1}}
\newcommand{\defnref}[1]{Definition~\ref{#1}}
\newcommand{\examref}[1]{Example~\ref{#1}}
\newcommand{\notatref}[1]{Notation~\ref{#1}}

\renewcommand{\eqref}[1]{(\ref{#1})}
\newcommand{\itemref}[1]{\ref{#1}}

\numberwithin{equation}{subsection}

\usepackage{tikz}
\usetikzlibrary{matrix}
\usepackage{tikz-cd}
\tikzset{
  commutative diagrams/.cd,
  arrow style=tikz,
  diagrams={>=latex}}
\makeatletter
\tikzset{
  column sep/.code=\def\pgfmatrixcolumnsep{\pgf@matrix@xscale*(#1)},
  row sep/.code   =\def\pgfmatrixrowsep{\pgf@matrix@yscale*(#1)},
  matrix xscale/.code=%
    \pgfmathsetmacro\pgf@matrix@xscale{\pgf@matrix@xscale*(#1)},
  matrix yscale/.code=%
    \pgfmathsetmacro\pgf@matrix@yscale{\pgf@matrix@yscale*(#1)},
  matrix scale/.style={/tikz/matrix xscale={#1},/tikz/matrix yscale={#1}}}
\def\pgf@matrix@xscale{1}
\def\pgf@matrix@yscale{1}
\makeatother

\usepackage{xspace}

\usepackage{xifthen}

\usepackage{xparse}

\usepackage{mathtools}

\newcommand{\nc}{\newcommand}
\nc{\renc}{\renewcommand}
\nc{\ssec}{\subsection}
\nc{\sssec}{\subsubsection}
\nc{\on}{\operatorname}
\nc{\term}[1]{#1\xspace}

\nc{\sA}{\ensuremath{\mathscr{A}}\xspace}
\nc{\sB}{\ensuremath{\mathscr{B}}\xspace}
\nc{\sC}{\ensuremath{\mathscr{C}}\xspace}
\nc{\sD}{\ensuremath{\mathscr{D}}\xspace}
\nc{\sE}{\ensuremath{\mathscr{E}}\xspace}
\nc{\sF}{\ensuremath{\mathscr{F}}\xspace}
\nc{\sG}{\ensuremath{\mathscr{G}}\xspace}
\nc{\sH}{\ensuremath{\mathscr{H}}\xspace}
\nc{\sI}{\ensuremath{\mathscr{I}}\xspace}
\nc{\sJ}{\ensuremath{\mathscr{J}}\xspace}
\nc{\sK}{\ensuremath{\mathscr{K}}\xspace}
\nc{\sL}{\ensuremath{\mathscr{L}}\xspace}
\nc{\sM}{\ensuremath{\mathscr{M}}\xspace}
\nc{\sN}{\ensuremath{\mathscr{N}}\xspace}
\nc{\sO}{\ensuremath{\mathscr{O}}\xspace}
\nc{\sP}{\ensuremath{\mathscr{P}}\xspace}
\nc{\sQ}{\ensuremath{\mathscr{Q}}\xspace}
\nc{\sR}{\ensuremath{\mathscr{R}}\xspace}
\nc{\sS}{\ensuremath{\mathscr{S}}\xspace}
\nc{\sT}{\ensuremath{\mathscr{T}}\xspace}
\nc{\sU}{\ensuremath{\mathscr{U}}\xspace}
\nc{\sV}{\ensuremath{\mathscr{V}}\xspace}
\nc{\sW}{\ensuremath{\mathscr{W}}\xspace}
\nc{\sX}{\ensuremath{\mathscr{X}}\xspace}
\nc{\sY}{\ensuremath{\mathscr{Y}}\xspace}
\nc{\sZ}{\ensuremath{\mathscr{Z}}\xspace}

\nc{\bA}{\ensuremath{\mathbf{A}}\xspace}
\nc{\bB}{\ensuremath{\mathbf{B}}\xspace}
\nc{\bC}{\ensuremath{\mathbf{C}}\xspace}
\nc{\bD}{\ensuremath{\mathbf{D}}\xspace}
\nc{\bE}{\ensuremath{\mathbf{E}}\xspace}
\nc{\bF}{\ensuremath{\mathbf{F}}\xspace}
\nc{\bG}{\ensuremath{\mathbf{G}}\xspace}
\nc{\bH}{\ensuremath{\mathbf{H}}\xspace}
\nc{\bI}{\ensuremath{\mathbf{I}}\xspace}
\nc{\bJ}{\ensuremath{\mathbf{J}}\xspace}
\nc{\bK}{\ensuremath{\mathbf{K}}\xspace}
\nc{\bL}{\ensuremath{\mathbf{L}}\xspace}
\nc{\bM}{\ensuremath{\mathbf{M}}\xspace}
\nc{\bN}{\ensuremath{\mathbf{N}}\xspace}
\nc{\bO}{\ensuremath{\mathbf{O}}\xspace}
\nc{\bP}{\ensuremath{\mathbf{P}}\xspace}
\nc{\bQ}{\ensuremath{\mathbf{Q}}\xspace}
\nc{\bR}{\ensuremath{\mathbf{R}}\xspace}
\nc{\bS}{\ensuremath{\mathbf{S}}\xspace}
\nc{\bT}{\ensuremath{\mathbf{T}}\xspace}
\nc{\bU}{\ensuremath{\mathbf{U}}\xspace}
\nc{\bV}{\ensuremath{\mathbf{V}}\xspace}
\nc{\bW}{\ensuremath{\mathbf{W}}\xspace}
\nc{\bX}{\ensuremath{\mathbf{X}}\xspace}
\nc{\bY}{\ensuremath{\mathbf{Y}}\xspace}
\nc{\bZ}{\ensuremath{\mathbf{Z}}\xspace}

\nc{\dA}{\ensuremath{\mathds{A}}\xspace}
\nc{\dB}{\ensuremath{\mathds{B}}\xspace}
\nc{\dC}{\ensuremath{\mathds{C}}\xspace}
\nc{\dD}{\ensuremath{\mathds{D}}\xspace}
\nc{\dE}{\ensuremath{\mathds{E}}\xspace}
\nc{\dF}{\ensuremath{\mathds{F}}\xspace}
\nc{\dG}{\ensuremath{\mathds{G}}\xspace}
\nc{\dH}{\ensuremath{\mathds{H}}\xspace}
\nc{\dI}{\ensuremath{\mathds{I}}\xspace}
\nc{\dJ}{\ensuremath{\mathds{J}}\xspace}
\nc{\dK}{\ensuremath{\mathds{K}}\xspace}
\nc{\dL}{\ensuremath{\mathds{L}}\xspace}
\nc{\dM}{\ensuremath{\mathds{M}}\xspace}
\nc{\dN}{\ensuremath{\mathds{N}}\xspace}
\nc{\dO}{\ensuremath{\mathds{O}}\xspace}
\nc{\dP}{\ensuremath{\mathds{P}}\xspace}
\nc{\dQ}{\ensuremath{\mathds{Q}}\xspace}
\nc{\dR}{\ensuremath{\mathds{R}}\xspace}
\nc{\dS}{\ensuremath{\mathds{S}}\xspace}
\nc{\dT}{\ensuremath{\mathds{T}}\xspace}
\nc{\dU}{\ensuremath{\mathds{U}}\xspace}
\nc{\dV}{\ensuremath{\mathds{V}}\xspace}
\nc{\dW}{\ensuremath{\mathds{W}}\xspace}
\nc{\dX}{\ensuremath{\mathds{X}}\xspace}
\nc{\dY}{\ensuremath{\mathds{Y}}\xspace}
\nc{\dZ}{\ensuremath{\mathds{Z}}\xspace}

\nc{\bbA}{\ensuremath{\mathbb{A}}\xspace}
\nc{\bbB}{\ensuremath{\mathbb{B}}\xspace}
\nc{\bbC}{\ensuremath{\mathbb{C}}\xspace}
\nc{\bbD}{\ensuremath{\mathbb{D}}\xspace}
\nc{\bbE}{\ensuremath{\mathbb{E}}\xspace}
\nc{\bbF}{\ensuremath{\mathbb{F}}\xspace}
\nc{\bbG}{\ensuremath{\mathbb{G}}\xspace}
\nc{\bbH}{\ensuremath{\mathbb{H}}\xspace}
\nc{\bbI}{\ensuremath{\mathbb{I}}\xspace}
\nc{\bbJ}{\ensuremath{\mathbb{J}}\xspace}
\nc{\bbK}{\ensuremath{\mathbb{K}}\xspace}
\nc{\bbL}{\ensuremath{\mathbb{L}}\xspace}
\nc{\bbM}{\ensuremath{\mathbb{M}}\xspace}
\nc{\bbN}{\ensuremath{\mathbb{N}}\xspace}
\nc{\bbO}{\ensuremath{\mathbb{O}}\xspace}
\nc{\bbP}{\ensuremath{\mathbb{P}}\xspace}
\nc{\bbQ}{\ensuremath{\mathbb{Q}}\xspace}
\nc{\bbR}{\ensuremath{\mathbb{R}}\xspace}
\nc{\bbS}{\ensuremath{\mathbb{S}}\xspace}
\nc{\bbT}{\ensuremath{\mathbb{T}}\xspace}
\nc{\bbU}{\ensuremath{\mathbb{U}}\xspace}
\nc{\bbV}{\ensuremath{\mathbb{V}}\xspace}
\nc{\bbW}{\ensuremath{\mathbb{W}}\xspace}
\nc{\bbX}{\ensuremath{\mathbb{X}}\xspace}
\nc{\bbY}{\ensuremath{\mathbb{Y}}\xspace}
\nc{\bbZ}{\ensuremath{\mathbb{Z}}\xspace}

\nc{\cA}{\ensuremath{\mathcal{A}}\xspace}
\nc{\cB}{\ensuremath{\mathcal{B}}\xspace}
\nc{\cC}{\ensuremath{\mathcal{C}}\xspace}
\nc{\cD}{\ensuremath{\mathcal{D}}\xspace}
\nc{\cE}{\ensuremath{\mathcal{E}}\xspace}
\nc{\cF}{\ensuremath{\mathcal{F}}\xspace}
\nc{\cG}{\ensuremath{\mathcal{G}}\xspace}
\nc{\cH}{\ensuremath{\mathcal{H}}\xspace}
\nc{\cI}{\ensuremath{\mathcal{I}}\xspace}
\nc{\cJ}{\ensuremath{\mathcal{J}}\xspace}
\nc{\cK}{\ensuremath{\mathcal{K}}\xspace}
\nc{\cL}{\ensuremath{\mathcal{L}}\xspace}
\nc{\cM}{\ensuremath{\mathcal{M}}\xspace}
\nc{\cN}{\ensuremath{\mathcal{N}}\xspace}
\nc{\cO}{\ensuremath{\mathcal{O}}\xspace}
\nc{\cP}{\ensuremath{\mathcal{P}}\xspace}
\nc{\cQ}{\ensuremath{\mathcal{Q}}\xspace}
\nc{\cR}{\ensuremath{\mathcal{R}}\xspace}
\nc{\cS}{\ensuremath{\mathcal{S}}\xspace}
\nc{\cT}{\ensuremath{\mathcal{T}}\xspace}
\nc{\cU}{\ensuremath{\mathcal{U}}\xspace}
\nc{\cV}{\ensuremath{\mathcal{V}}\xspace}
\nc{\cW}{\ensuremath{\mathcal{W}}\xspace}
\nc{\cX}{\ensuremath{\mathcal{X}}\xspace}
\nc{\cY}{\ensuremath{\mathcal{Y}}\xspace}
\nc{\cZ}{\ensuremath{\mathcal{Z}}\xspace}

\nc{\mrm}[1]{\ensuremath{\mathrm{#1}}\xspace}
\nc{\mbf}[1]{\ensuremath{\mathbf{#1}}\xspace}
\nc{\mcal}[1]{\ensuremath{\mathcal{#1}}\xspace}
\nc{\mit}[1]{\ensuremath{\mathit{#1}}\xspace}
\nc{\msc}[1]{\ensuremath{\mathscr{#1}}\xspace}

\let\sectsign\S
\let\S\relax

\nc{\sub}{\subseteq}
\nc{\too}{\longrightarrow}
\nc{\hook}{\hookrightarrow}
\nc{\hooklongrightarrow}{\lhook\joinrel\longrightarrow}
\nc{\hooklong}{\hooklongrightarrow}
\nc{\twoheadlongrightarrow}{\relbar\joinrel\twoheadrightarrow}
\nc{\longrightleftarrows}{\ \raisebox{0.3ex}{\(\mathrel{\substack{\xrightarrow{\rule{1em}{0em}} \\[-1ex] \xleftarrow{\rule{1em}{0em}}}}\)}\ }

\renc{\ge}{\geqslant}
\renc{\le}{\leqslant}

\nc{\shiso}{\approx}
\nc{\isoto}{\stackrel{\sim}{\smash{\longrightarrow}\rule{0pt}{0.4ex}}}
\nc{\isofrom}{\xleftarrow{\sim}}
\usepackage{stackengine}
\nc{\altxrightarrow}[2][0pt]{\mathrel{\ensurestackMath{\stackengine%
  {\dimexpr#1-5.5pt}{\xrightarrow{\phantom{#2}}}{\scriptstyle\!#2\,}%
  {O}{c}{F}{F}{S}}}}
\nc{\isotoo}{\altxrightarrow{~\sim~}}

\nc{\id}{\mathrm{id}}

\DeclareMathOperator{\Hom}{\mathrm{Hom}}
\nc{\uHom}{\underline{\smash{\Hom}}}
\DeclareMathOperator{\Maps}{\mathrm{Maps}}

\DeclareMathOperator{\End}{\mathrm{End}}

\nc{\uEnd}{\underline{\smash{\End}}}

\nc{\colim}{\varinjlim}
\renc{\lim}{\varprojlim}
\nc{\Cofib}{\on{Cofib}}
\nc{\Fib}{\on{Fib}}
\nc{\initial}{\varnothing}
\nc{\op}{\mathrm{op}}

\DeclareMathOperator*{\fibprod}{\times}

\nc{\Spc}{{\mrm{Spc}}}
\nc{\Spt}{{\mrm{Spt}}}
\nc{\pt}{{\mrm{pt}}}
\nc{\Stab}{\mrm{Stab}}
\nc{\Fun}{\mrm{Fun}}
\nc{\h}{\mrm{h}}
\renc{\L}{\mbf{L}}

\nc{\CAlg}{\mrm{CAlg}}
\nc{\cn}{{\mrm{cn}}}

\nc{\Spec}{\on{Spec}}
\nc{\A}{\bA}
\renc{\P}{\bP}
\nc{\Sm}{\mrm{Sm}}
\nc{\Aff}{\mrm{Aff}}
\nc{\cl}{{\mrm{cl}}}
\nc{\Zar}{\mathrm{Zar}}
\nc{\Nis}{\mathrm{Nis}}
\nc{\sm}{\mrm{sm}}

\nc{\Qcoh}{\on{Qcoh}}
\nc{\Perf}{\on{Perf}}

\nc{\MotSpc}{{\mbf{H}}}
\nc{\LNis}{\on{\mrm{L}_\Nis}}
\nc{\LZar}{\on{\mrm{L}_\Zar}}
\nc{\htp}{{\A^1}}
\nc{\Lhtp}{\on{\mrm{L}_\htp}}

\nc{\hspc}[2][]{\h_{#1}(#2)}

\nc{\nil}{{\mrm{nil}}}

\nc{\SPC}{{\underline{\smash{\mrm{Spc}}}}}
\nc{\MOTSPC}{{\underline{\smash{\mbf{H}}}}}

\nc{\AffCl}{\mrm{AffCl}}
\nc{\SmCl}{\mrm{SmCl}}
\nc{\htpcl}{{\A^1_\cl}}
\nc{\Lhtpcl}{\mrm{L}_{\htpcl}}

\nc{\fibsm}{\mathrm{fibsm}}
\nc{\htpflat}{{\A^{1,\flat}}}
\nc{\Lhtpflat}{\mrm{L}_{\A^{1,\flat}}}
\nc{\Lmotflat}{\L^{\flat}}

\nc{\Einfty}{{\sE_\infty}}
\nc{\bDelta}{\mbf{\Delta}}
\nc{\Tot}{\on{Tot}}
\nc{\Cech}{\textnormal{\v{C}}}
\nc{\T}{\bT}

\nc{\B}{\mrm{B}}
\nc{\K}{\mrm{K}}
\nc{\KB}{\mrm{K}^\B}
\nc{\KH}{\mrm{KH}}
\nc{\Mod}{\mrm{Mod}}
\nc{\perf}{\mrm{perf}}
\nc{\proj}{\mrm{proj}}
\nc{\free}{\mrm{free}}
\nc{\G}{\mbf{G}}
\nc{\GL}{\mrm{GL}}
\nc{\BGL}{\mrm{BGL}}
\nc{\gp}{\mrm{gp}}
\nc{\Q}{\on{Q}}

\nc{\Ering}{\term{$\Einfty$-ring}}
\nc{\Erings}{\term{$\Einfty$-rings}}

\nc{\cnEring}{\term{connective $\Einfty$-ring}}
\nc{\cnErings}{\term{connective $\Einfty$-rings}}

\nc{\inftyCat}{\term{$\infty$-category}}
\nc{\inftyCats}{\term{$\infty$-categories}}

\nc{\inftyGrpd}{\term{$\infty$-groupoid}}
\nc{\inftyGrpds}{\term{$\infty$-groupoids}}

\date{\today}
\title{\mbox{$\mathbf{A}^1$-homotopy invariance in spectral algebraic geometry}}
\author{Denis-Charles Cisinski}
\address{Fakultät für Mathematik\\
Universität Regensburg\\
93040 Regensburg\\
Germany}
\email{\href{mailto:denis-charles.cisinski@mathematik.uni-regensburg.de}{denis-charles.cisinski@mathematik.uni-regensburg.de}}
\urladdr{\url{http://www.mathematik.uni-regensburg.de/cisinski/}}

\author{Adeel A. Khan}
\address{Institut des Hautes Études Scientifiques\\
35 route de Chartres\\
91440 Bures-sur-Yvette\\
France}
\email{\href{mailto:khan@ihes.fr}{khan@ihes.fr}}
\urladdr{\url{https://www.preschema.com}}

\begin{document}

\begin{abstract}
  We study two different flavours of $\A^1$-homotopy theory in the setting of spectral algebraic geometry, and compare them to classical $\A^1$-homotopy theory.
  As an application we show that the spectral analogue of Weibel's homotopy invariant K-theory collapses to the classical theory.
  Along the way we give a new construction of nonconnective algebraic K-theory of stable \inftyCats via a generalization of the Bass--Thomason--Trobaugh construction.
\end{abstract}

\maketitle

\parskip 0pt
\tableofcontents

\renc{\baselinestretch}{1.1}
\setlength{\parindent}{0em}
\parskip 1em


\section{Introduction}
\label{sec:intro}

This paper establishes in a systematic way why fundamental invariants from derived geometry, such as Serre's Tor formula, or virtual fundamental classes, have a natural interpretation in homotopy invariant (co)homology theories.
In fact, we provide an explicit way to interpret Lurie's spectral geometry into Voevodsky's motivic homotopy theory as follows.

Let $R$ be a commutative ring or connective $\Einfty$-ring spectrum, and let $R[T]$ denote the polynomial algebra over $R$ in one variable $T$.
A peculiarity of the world of $\Einfty$-ring spectra is that the polynomial algebra $R[T]$ is \emph{free} as an $\Einfty$-$R$-algebra only when $R$ is of characteristic zero (an $\Einfty$-$\bQ$-algebra).
If we write $R\{T\}$ for the free $\Einfty$-$R$-algebra on one generator $T$ (in degree zero), then in general there is only a comparison homomorphism $R\{T\} \to R[T]$.
The $\Einfty$-ring $R\{T\}$ is \emph{smooth} over $R$ (its cotangent complex $L_{R\{T\}/R}$ is free) but not flat, while the $\Einfty$-ring $R[T]$ is usually not smooth but instead \emph{fibre-smooth}: that is, it is flat over $R$, and $\pi_0(R[T]) \simeq \pi_0(R)[T]$ is smooth over $\pi_0(R)$ in the sense of ordinary commutative algebra.
Let $\CAlg^\sm_R$ denote the \inftyCat of smooth $\Einfty$-algebras over $R$, and $\CAlg^\fibsm_{R}$ denote the \inftyCat of fibre-smooth $\Einfty$-algebras over $R$.
Our first main result reads as follows:

\begin{thmX}\label{thm:intro/comp1}
  Let $R$ be a \cnEring.
  Consider the following \inftyCats:
  \begin{thmlist}
    \item
    The \inftyCat $\MotSpc(R)$ of Nisnevich sheaves of spaces $\sF : \CAlg^\sm_R \to \Spc$ for which the canonical map $\sF(A) \to \sF(A\{T\})$ is invertible for every $A \in \CAlg^\sm_R$.

    \item
    The \inftyCat $\MotSpc^\flat(R)$ of Nisnevich sheaves of spaces $\sF : \CAlg^\fibsm_R \to \Spc$ for which the canonical map $\sF(A) \to \sF(A[T])$ is invertible for every $A \in \CAlg^\fibsm_R$.
    
    \item
    The \inftyCat $\MotSpc^\cl(\pi_0(R))$ of Nisnevich sheaves of spaces $\sF : \CAlg^\fibsm_{\pi_0(R)} \to \Spc$ for which the canonical map $\sF(A) \to \sF(A[T])$ is invertible for every $A \in \CAlg^{\fibsm}_{\pi_0(R)}$.
  \end{thmlist}
  Then (i) and (iii) are equivalent.
  If $R$ is an $\Einfty$-$\bZ$-algebra, then all three are equivalent.
\end{thmX}

See Theorems~\ref{thm:comp} and \ref{thm:comp2}.
The result is nontrivial even for ordinary commutative rings $R$ (viewed as discrete \Erings), in which case it asserts an equivalence $\MotSpc(R) \simeq \MotSpc^\cl(R)$.
Note that $\CAlg^\fibsm_{R}$ coincides with the category of commutative rings that are smooth over $R$ in the sense of ordinary commutative algebra, so that $\MotSpc^\cl(R)$ coincides with the usual $\A^1$-homotopy category considered by Morel and Voevodsky \cite{MorelVoevodsky}.
This equivalence was only known in characteristic zero (see Proposition~2.4.6 and Warning~2.4.7 in \cite{KhanLocalization}).
For simplicial commutative rings (regarded as $\Einfty$-$\bZ$-algebras), the equivalence (ii) $\Leftrightarrow$ (iii) was proven in the second author's thesis \cite{KhanThesis}.

Our second subject of discussion is a variant of Weibel's homotopy invariant K-theory for \cnErings.
Let $\KH^\cl$ denote the classical variant \cite{WeibelKH}, defined by starting with nonconnective algebraic K-theory and forcing it to become $\A^1$-homotopy invariant in the sense that the canonical map
\[ \KH^\cl(R) \to \KH^\cl(R[T]) \]
is invertible for every commutative ring $R$.
We define an analogous construction on the \inftyCat of \cnErings,
\begin{equation*}
  R \mapsto \KH(R),
\end{equation*}
by forcing $\A^1$-homotopy invariance in the sense that $\KH(R\{T\}) \to \KH(R)$ is invertible for any \cnEring $R$.
We then have the following comparison, a K-theoretic incarnation of \thmref{thm:intro/comp1}:

\begin{thmX}\label{thm:intro/KH}
  For every \cnEring $R$, there is a canonical isomorphism of spectra
  \begin{equation*}
    \KH(R) \simeq \KH^\cl(\pi_0(R)).
  \end{equation*}
\end{thmX}

See \corref{cor:KH global sections}.
For connective $\Einfty$-$\bZ$-algebras one may adapt the proof, using the equivalence (i) $\Leftrightarrow$ (iii) in \thmref{thm:intro/comp1}, to derive the same result for the variant $\KH^\flat$ constructed by imposing invertibility of the maps $\KH^\flat(R) \to \KH^\flat(R[T])$.
This result was communicated to us by B.~Antieau and D.~Gepner in 2015 in the generality of connective $\sE_1$-rings, and has recently been recorded in the case of connective $\sE_1$-$\bZ$-algebras by Land--Tamme \cite[Prop.~3.14]{LandTamme}.

An important ingredient in the proof is the observation that nonconnective K-theory of stable \inftyCats, defined as in \cite{BlumbergGepnerTabuada} following Schlichting, can also be described by a variant of the Bass--Thomason--Trobaugh construction \cite[Sect.~6]{ThomasonTrobaugh} defined over the sphere spectrum:

\begin{thmX}\label{thm:intro/Bass}
  There is an isomorphism \[ \mathbb{K} \simeq \K^\B \] of spectrum-valued functors on the \inftyCat of small stable \inftyCats, where $\K$ is algebraic K-theory, $\mathbb{K}$ is nonconnective algebraic K-theory, and $(-)^\B$ denotes a generalization of the Bass--Thomason--Trobaugh construction (\constrref{constr:Bass}).
\end{thmX}

See \examref{exam:K^B}.
In fact, we show in \thmref{thm:deloop localizing invariant} that the construction $(-)^\B$ defines an equivalence between the \inftyCat of connective spectrum-valued localizing invariants\footnote{%
  A connective spectrum-valued invariant $E$ is localizing if, for any short exact sequence of stable \inftyCats $\bA' \to \bA \to \bA''$, $E(\bA')$ is identified with the connective cover of the homotopy fibre of $E(\bA) \to E(\bA'')$.
  See \defnref{defn:localizing connective}.
} and the \inftyCat of spectrum-valued localizing invariants.
This also implies for example that all operations on connective K-theory deloop to $\mathbb{K}$.
The result is inspired by closely related work of Robalo in his framework of noncommutative motivic homotopy theory \cite{RobaloBridge} which gives similar results in the more restrictive setting of dg-categories.

\ssec*{Outline}

  In the body of the paper, we use the language of spectral algebraic geometry \cite{SAG-20180204}.
  Given a spectral affine\footnote{%
    To simplify the exposition we usually only discuss the affine case.
    However, all our results extend to spectral schemes and algebraic spaces by descent: see Corollaries~\ref{cor:comp asp}, \ref{cor:comp2 asp}, and \ref{cor:KH global sections}.
  } scheme $S$, we may define $\MotSpc(S)$ as the \inftyCat of $\A^1$-invariant Nisnevich sheaves on the site $\Sm_{/S}$ of smooth spectral affine schemes over $S$.
  Then $\MotSpc(\Spec(R))$ is equivalent to $\MotSpc(R)$ as defined above, and also to the construction given in \cite{KhanLocalization}, see Corollary~2.4.5 of \emph{op. cit}.
  Similarly we have the variant $\MotSpc^\flat(S)$ defined as the \inftyCat of $\A^{1,\flat}$-invariant Nisnevich sheaves on the site $\Sm^\flat_{/S}$ of fibre-smooth spectral affine schemes over $S$, where $\A^{1,\flat} = \Spec(\bS[T])$.

  The proof of the equivalence (i) $\Leftrightarrow$ (iii) in \thmref{thm:intro/comp1} is given in \secref{sec:nil}.
  The key idea is to re-interpret the derived nil-invariance result of \cite[Thm.~A]{KhanLocalization} as a sort of descent statement with respect to the ``nil topology'' whose coverings are morphisms of the form $X_\cl \to X$, where $X_\cl$ denotes the classical truncation of the spectral scheme $X$.
  Since the site $\Sm_{/S}$ is typically not closed under the operation $X \mapsto X_\cl$, making sense of this idea requires us to enlargen our site\footnote{%
    The same thing happens in classical algebraic geometry with the cdh topology; see \cite{KhanCdh}.
  }.
  The first few subsections (\ref{ssec:nil/Asp-fibred}, \ref{ssec:nil/extension}, \ref{ssec:nil/functoriality}) develop some generalities related to $\A^1$-homotopy theory on various \emph{admissible} sites (\defnref{defn:adm}) and how variation of site interplays with the basic operations such as inverse/direct image, product, and internal hom.
  For example, any \emph{narrow} subcategory $\sA_{/S} \sub \Aff_{/S}$ as in \defnref{defn:narrow} gives rise to the same $\A^1$-homotopy theory as the smooth site $\Sm_{/S}$ (\examref{exam:narrow=Sm-fibred}).
  The key result here is \propref{prop:MotSpc vs MOTSPC} which implies that any $\A^1$-invariant Nisnevich sheaf defined on a narrow subcategory $\sA_{/S}$ admits a canonical extension $\sF^+$ to an $\A^1$-invariant Nisnevich sheaf on any \emph{broad} site $\sB_{/S}$ (\defnref{defn:broad}) that contains $\sA_{/S}$.
  Broad sites are in particular closed under classical truncation.
  
  In \ssecref{ssec:nil/classical} we construct a comparison functor to the classical motivic homotopy category, which we show in \ssecref{ssec:nil/nil-localization} is a left Bousfield localization (the \emph{nil-localization}) if we work over a broad site (\thmref{thm:nil-localization}).
  In \ssecref{ssec:nil/nil-local} we formulate the nil-descent result alluded to above (\thmref{thm:nil descent}).
  Finally, we put everything together in \ssecref{ssec:nil/comparison} to show that nil-localization gives an equivalence
  \[
    \MotSpc^\cl(S_\cl) \isotoo \MotSpc(S)
  \]
  when we restrict to our narrow subcategory $\sA_{/S}$.
  The last subsection (\ssecref{ssec:nil/V-linear}) extends to the result to presheaves with values in general presentable \inftyCats $\bV$ (e.g. presheaves of spectra or presheaves of chain complexes).

  The second part of \thmref{thm:intro/comp1} is proven in \secref{sec:smflat}.
  The idea of enlarging sites again plays an important role here.
  Working on a site large enough that it contains both the spectral affine line $\A^1$ and the flat affine line $\A^{1,\flat}$ allows us to exploit the canonical morphism $\varepsilon : \A^{1,\flat} \to \A^1$ which is a morphism of interval objects (\ref{rem:A^1_flat to A^1 intervals}).
  The key input in the comparison is the fact that $\varepsilon$ is a ``universal'' $\A^{1,\flat}$-equivalence, as long as we work over $\Spec(\bZ)$ (\lemref{lem:epsilon A^1_flat-equiv}).
  The proof of the comparison is given in \ssecref{ssec:smflat/comp}.

  \secref{sec:bass} is independent of the first three sections and discusses localizing invariants of stable \inftyCats in the sense of \cite{BlumbergGepnerTabuada}\footnote{%
    Note that, unlike \cite{BlumbergGepnerTabuada}, we do not require localizing invariants to preserve filtered colimits.
  }.
  The main result, \thmref{thm:deloop localizing invariant}, asserts that every connective spectrum-valued localizing invariant admits a unique delooping to a spectrum-valued localizing invariant.
  This is proven by generalizing the Bass construction to stable \inftyCats over the sphere spectrum (or any \cnEring).

  \thmref{thm:intro/KH} is proven in \secref{sec:kh} by showing that the equivalence $\MotSpc(R) \simeq \MotSpc^\cl(\pi_0(R))$ of \thmref{thm:intro/comp1} or (rather its generalization to sheaves of spectra) sends $\KH$ to $\KH^\cl$.
  Unstably this boils down to representability results for the infinite loop spaces, and the stable result is deduced via Bott periodicity.

\ssec*{Notation and conventions}

  We will use the language of \inftyCats freely throughout the text.
  Our main references are \cite{HTT,HA-20170918}.
  The \inftyCat of spaces and spectra will be denoted by $\Spc$ and $\Spt$, respectively, and a morphism in an \inftyCat will be called an \emph{isomorphism} if it is invertible (= an \emph{equivalence} in the language of \cite{HTT}).
  We also use the language of spectral algebraic geometry \cite{SAG-20180204}.
  Given a spectral affine scheme $S = \Spec(R)$, we write $\Aff_{/S}$ for the \inftyCat of spectral affine schemes over $S$, which is equivalent to the opposite of the \inftyCat of connective $\Einfty$-$R$-algebras.

\ssec*{Acknowledgments}

  We would like to thank Benjamin Antieau and David Gepner for their encouragement and interest in this paper.
  The second-named author thanks Institut Mittag--Leffler for hospitality while the first draft of this paper was written.
  He also acknowledges partial support from SFB 1085 Higher Invariants, Universit\"at Regensburg, and from the Simons Collaboration on Homological Mirror Symmetry.


\section{Comparison with classical motivic homotopy theory}
\label{sec:nil}


\ssec{Fibred spaces}
\label{ssec:nil/Asp-fibred}

For this subsection, we fix an affine spectral scheme $S$ and write $\Aff_{/S}$ for the \inftyCat of affine spectral schemes over $S$.

\begin{defn}\label{defn:smooth}
  A morphism of affine spectral schemes $X \to S$ is called \emph{smooth} (resp. \emph{étale}) if it is of finite presentation and the relative cotangent complex $\sL_{X/S}$ is a locally free $\sO_X$-module of finite rank (resp. is zero).
\end{defn}

\begin{rem}
  From \cite[Prop.~11.2.2.1]{SAG-20180204} it follows that a morphism of affine spectral schemes $\Spec(B) \to \Spec(A)$ is smooth if and only if $A \to B$ is differentially smooth in the sense of \cite[Def.~11.2.2.2]{SAG-20180204}.
\end{rem}

\begin{exam}\label{exam:A^n}
  Let $\bS$ denote the sphere spectrum.
  For every integer $n\ge 0$, we write $\bS\{T_1,\ldots,T_n\}$ for the free $\Einfty$-algebra on $n$ generators $T_i$ (in degree zero).
  We let $\A^n$ denote the affine spectral scheme $\Spec(\bS\{T_1,\ldots,T_n\})$ and refer to it as $n$-dimensional \emph{spectral affine space} (over the sphere spectrum).
  The morphism $\A^n \to \Spec(\bS)$ is smooth, and we have a canonical isomorphism $(\A^n)_\cl \simeq \A^n_\cl$, where $\A^n_\cl$ denotes the classical affine space over $\Spec(\bZ)$.
\end{exam}

\begin{rem}\label{rem:smooth factor}
  If $X \in \Aff_{/S}$ is smooth over $S$, then Zariski-locally on $X$ there exists an étale $S$-morphism $X \to S \times \A^n$ for some $n\ge 0$.
  This follows from \cite[Prop.~11.2.2.1]{SAG-20180204}.
\end{rem}

\begin{defn}[Nisnevich excision]\label{defn:Nisnevich}\leavevmode
  \begin{defnlist}
    \item
    A \emph{Nisnevich square} over $X \in \Aff_{/S}$ is a cartesian square of affine spectral schemes
      \begin{equation} \label{eq:Nisnevich square}
        \begin{tikzcd}
          W \arrow{r}\arrow{d}
            & V \arrow{d}{p}
          \\
          U \arrow{r}{j}
            & X
        \end{tikzcd}
      \end{equation}
    where $j$ is an open immersion, $p$ is \'etale, and there exists a closed immersion $Z \hook X$ complementary to $j$ such that the induced morphism $p^{-1}(Z) \to Z$ is invertible.

    \item
    We say that a presheaf of spaces $\sF$ on $\Aff_{/S}$ satisfies \emph{Nisnevich excision} if it is reduced, i.e. the space $\Gamma(\initial, \sF)$ is contractible, and for any $X \in \Aff_{/S}$ and any Nisnevich square over $X$ of the form \eqref{eq:Nisnevich square}, the induced square of spaces
      \begin{equation*}
        \begin{tikzcd}
          \Gamma(X, \sF) \arrow{r}{j^*}\arrow{d}{p^*}
            & \Gamma(U, \sF) \arrow{d}
          \\
          \Gamma(V, \sF) \arrow{r}
            & \Gamma(W, \sF)
      \end{tikzcd}
      \end{equation*}
    is cartesian.
  \end{defnlist}
\end{defn}

\begin{defn}[$\A^1$-invariance]\label{defn:A1}
  Let $\sF$ be a presheaf of spaces on $\Aff_{/S}$.
  We say that $\sF$ satisfies \emph{$\A^1$-homotopy invariance} if for every $X \in \Aff_{/S}$, the canonical map of spaces
  \begin{equation*}
    p^* : \Gamma(X, \sF) \to \Gamma(X\times\A^1, \sF)
  \end{equation*}
  is invertible, where $p : X\times\A^1 \to X$ is the projection of the spectral affine line over $X$.
\end{defn}

We will need to consider presheaves defined on smaller subcategories of $\Aff_{/S}$.
The following definition identifies the minimal conditions under which Definitions~\ref{defn:Nisnevich} and \ref{defn:A1} make sense.

\begin{defn}\label{defn:adm}
  We say that a full subcategory $\sC_{/S} \sub \Aff_{/S}$ is \emph{admissible} if it is essentially small and satisfies the following conditions:
  \begin{defnlist}[series=adm]
    \item\label{item:adm1}
    The affine spectral scheme $S$ (viewed over $S$ via the identity) belongs to $\sC_{/S}$.
    \item\label{item:adm2}
    If $X$ belongs to $\sC_{/S}$ and $Y$ is étale over $X$, then $Y$ belongs to $\sC_{/S}$.
    \item\label{item:adm3}
    If $X$ belongs to $\sC_{/S}$, then $X \times \A^n$ belongs to $\sC_{/S}$ for every $n\ge 0$.
  \end{defnlist}
\end{defn}

\begin{exam}\label{exam:Sm}
  The full subcategory $\Sm_{/S} \sub \Aff_{/S}$ of \emph{smooth} affine spectral $S$-schemes is admissible.
  This follows from the fact that étale morphisms are smooth, the morphism $\A^n \to \Spec(\bS)$ is smooth for every $n\ge 0$, and the class of smooth morphisms is stable under composition and base change.
\end{exam}

\begin{exam}\label{exam:minimal adm}
  Let $\sA^0_{/S} \sub \Aff_{/S}$ denote the full subcategory spanned by $X \in \Aff_{/S}$ which admit an étale morphism
  \[ X \to S \times \A^n \]
  over $S$.
  Then $\sA^0_{/S}$ is admissible.
  For the third condition, note that if $X$ admits an étale $S$-morphism to $S \times \A^n$, then $X \times \A^m$ admits (for every $m$) an étale $S$-morphism
  \[
    X \times \A^m
    \to S \times \A^n \times \A^m
    \xrightarrow{\mrm{pr}} S \times \A^n,
  \]
  and hence also belongs to $\sA^0_{/S}$.
  Note that $\sA^0_{/S}$ is in fact the minimal admissible subcategory of $\Aff_{/S}$.
\end{exam}

\begin{defn}\label{defn:narrow}
  Let $\sA_{/S} \sub \Aff_{/S}$ be an admissible subcategory.
  We say that $\sA_{/S}$ is \emph{narrow} if it is contained in the full subcategory $\Sm_{/S}$.
\end{defn}

\begin{defn}\label{defn:broad}
  Let $\sA_{/S} \sub \Aff_{/S}$ be an admissible subcategory.
  We say that $\sA_{/S}$ is \emph{broad} if it satisfies the following further condition:
  \begin{defnlist}[resume=adm]
    \item\label{item:adm4}
    For every $X \in \sA_{/S}$, the classical truncation $X_\cl$ also belongs to $\sA_{/S}$.
  \end{defnlist}
  Note that there is a minimal broad subcategory of $\Aff_{/S}$, which is the closure of \examref{exam:minimal adm} under the operations \ref{item:adm2}, \ref{item:adm3}, and \ref{item:adm4} (constructed by transfinite iteration).
\end{defn}

\begin{rem}
  Note that any broad subcategory contains the minimal narrow subcategory $\sA^0_{/S}$ (\examref{exam:minimal adm}).
  Note also that, as long as $S$ is not discrete, no admissible subcategory is both narrow and broad, since the morphism $S_\cl \to S$ is smooth if and only if it is an isomorphism.
\end{rem}

\begin{defn}\label{defn:fibred space}
  Let $\sC_{/S} \sub \Aff_{/S}$ be a full subcategory.
  A \emph{$\sC$-fibred space} over $S$ is a presheaf of spaces on $\sC_{/S}$.
\end{defn}

\begin{defn}\label{defn:A-fibred motivic space}
  Let $\sC_{/S} \sub \Aff_{/S}$ be an admissible subcategory.
  We say that a $\sC$-fibred space $\sF$ satisfies \emph{Nisnevich excision} and \emph{$\A^1$-homotopy invariance} if it satisfies the conditions of Definitions~\ref{defn:Nisnevich} and \ref{defn:A1}, respectively (imposed only on objects $X\in\sC_{/S}$).
  We say that $\sF$ is a \emph{$\sC$-fibred motivic space} over $S$ if it is both Nisnevich excisive and $\A^1$-homotopy invariant.
  
  We denote by
  \[ \Spc(\sC_{/S})\quad \text{and}~ \MotSpc(\sC_{/S}) \]
  the \inftyCat of $\sC$-fibred spaces over $S$ and its full subcategory of motivic objects.
\end{defn}

\begin{exam}\label{exam:MotSpc(S)}
  In case of the admissible subcategory $\Sm_{/S} \sub \Aff_{/S}$ (\examref{exam:Sm}), we write
  \[
    \Spc(S) := \Spc(\Sm_{/S}),
    \quad \MotSpc(S) := \MotSpc(\Sm_{/S}).
  \]
  With this definition we have $\MotSpc(\Spec(R)) \simeq \MotSpc(R)$ for any \cnEring $R$, where the right-hand side is as defined in \thmref{thm:intro/comp1}.
\end{exam}

\begin{rem}\label{rem:SPC localization basics}
  For any admissible subcategory $\sC_{/S} \sub \Aff_{/S}$, the full subcategories of Nisnevich-excisive, $\A^1$-invariant, and motivic $\sC$-fibred spaces are each left Bousfield localizations of the \inftyCat of $\sC$-fibred spaces.
  The following assertions are proven in the same way as their analogues for $\Sm$-fibred spaces (cf. \cite[Sect.~2]{KhanLocalization}):
  \begin{defnlist}
    \item
    The Nisnevich localization functor $\sF \mapsto \LNis(\sF)$ is exact (follows from \cite[Thm.~2.2.7]{KhanLocalization}).

    \item
    The $\A^1$-localization functor $\sF \mapsto \Lhtp(\sF)$ admits the following description (see \cite[Rem.~2.3.5]{KhanLocalization}, \cite[Prop.~3.4]{HoyoisEquivariant}): for every $\sC$-fibred space $\sF$, the space of sections over any $X \in \sC_{/S}$ is computed by a sifted colimit:
      \begin{equation} \label{eq:formula for Lhtp}
        \Gamma\big(X, \Lhtp(\sF)\big) \simeq \colim \Gamma(\A^n \times X, \sF),
      \end{equation}
    indexed by the opposite of the full subcategory $\bA_X \sub \Aff_{/X}$ whose objects are spectral affine spaces $X\times\A^n$ ($n\ge 0$).

    \item
    The motivic localization functor $\sF \mapsto \L(\sF)$ can be computed as the transfinite composite
      \begin{equation} \label{eq:description of Lmot}
        \L(\sF) \simeq \colim_{n \ge 0} (\Lhtp \circ \LNis)^{\circ n}(\sF),
      \end{equation}
    for any $\sF \in \Spc(\sC_{/S})$ (cf. \cite[Rem.~2.4.3]{KhanLocalization}).
    Moreover, $\MotSpc(\sC_{/S})$ has universality of colimits.

    \item
    The \inftyCat $\MotSpc(\sC_{/S})$ of $\sC$-fibred motivic spaces is generated under sifted colimits by objects of the form $\L \h_S(X)$, where $\h_S(X)$ is the presheaf on $\sC_{/S}$ represented by $X \in \sC_{/S}$ (cf. \cite[Prop.~2.4.4]{KhanLocalization}).
  \end{defnlist}
\end{rem}


\ssec{Extension of fibred spaces}
\label{ssec:nil/extension}

As in \ssecref{ssec:nil/Asp-fibred}, we fix an affine spectral scheme $S$.
We also fix the following data:

\begin{notat}\label{notat:inclusion}
  Fix an inclusion $\sC_{/S} \sub \sD_{/S}$ of admissible subcategories of $\Aff_{/S}$.
  We consider the \inftyCats
  \[ \Spc(\sC_{/S})\quad \text{and}~ \MotSpc(\sC_{/S}), \]
  \[ \Spc(\sD_{/S})\quad \text{and}~ \MotSpc(\sD_{/S}), \]
  as in \defnref{defn:A-fibred motivic space}.
\end{notat}

In this subsection we show that there are fully faithful embeddings
\[
  \Spc(\sC_{/S}) \hook \Spc(\sD_{/S}),
  \quad \MotSpc(\sC_{/S}) \hook \MotSpc(\sD_{/S}).
\]

\begin{notat}\label{notat:Sm -> Sm+}
  Let $\iota : \sC_{/S} \hook \sD_{/S}$ denote the inclusion functor.
  Restriction along $\iota$ defines a functor $\iota^* : \Spc(\sD_{/S}) \to \Spc(\sC_{/S})$, whose left adjoint $\iota_! : \Spc(\sC_{/S}) \to \Spc(\sD_{/S})$ is given by left Kan extension of $\iota$.
  The latter is uniquely characterized by the property of commutativity with colimits, and the identity $\iota_!\hspc[S]{X} \simeq \hspc[S]{X}$ for $X \in \sC_{/S}$.
  In particular, $\iota_!$ is fully faithful with essential image generated under colimits by objects of the form $\hspc[S]{X}$, with $X \in \sC_{/S}$.
  Similarly, $\iota^*$ also admits a fully faithful right adjoint $\iota_*$ given by right Kan extending $\iota$.
\end{notat}

\begin{prop}\label{prop:MotSpc vs MOTSPC}
In the notation of \ref{notat:Sm -> Sm+}, the assignment $\sF \mapsto \L\iota_!(\sF)$ induces a fully faithful functor of \inftyCats
  \begin{equation*}
    \L\iota_! : \MotSpc(\sC_{/S}) \to \MotSpc(\sD_{/S}),
  \end{equation*}
whose essential image is generated under sifted colimits by objects of the form $\L \h_S(X)$, where $X$ belongs to $\sC_{/S}$.
\end{prop}

\begin{exam}\label{exam:narrow=Sm-fibred}
  Suppose $\sA_{/S} \sub \Aff_{/S}$ is a \emph{narrow} subcategory, and consider the inclusion $\iota : \sA_{/S} \hook \Sm_{/S}$.
  Then the fully faithful embedding
  \[ \L\iota_! : \MotSpc(\sA_{/S}) \hook \MotSpc(\Sm_{/S}) = \MotSpc(S) \]
  is an \emph{equivalence}.
  Indeed, by \cite[Prop.~2.4.4]{KhanLocalization}, $\MotSpc(S)$ is generated under sifted colimits by objects of the form $\L \h_S(X)$, where $X$ belongs to the minimal admissible subcategory $\sA^0_{/S}$ (\examref{exam:minimal adm}), and hence also to $\sA_{/S}$.
\end{exam}

We will deduce \propref{prop:MotSpc vs MOTSPC} from the following lemma:

\begin{lem}\label{lem:iota and local equivs}
  The functors $\iota_! : \Spc(\sC_{/S}) \to \Spc(\sD_{/S})$ and $\iota^* : \Spc(\sD_{/S}) \to \Spc(\sC_{/S})$ preserve Nisnevich-local and $\A^1$-local equivalences.
\end{lem}

\begin{proof}
Since $\iota$ preserves Nisnevich squares and $\A^1$-projections, it follows that $\iota_!$ preserves Nisnevich-local and $\A^1$-local equivalences.
The definition of admissibility implies that $\iota$ is also cocontinuous with respect to the Nisnevich topology, so it follows that $\iota^*$ preserves Nisnevich-local equivalences (see \cite[Exp. III, Prop.~2.2]{SGA4} or \cite[Def.~3.1.5]{KhanLocalization}).

For $\A^1$-local equivalences it will suffice to show that, for any $X\in\sD_{/S}$, the canonical morphism
  \begin{equation*}
    \iota^*\h_S(X \times \A^1) \to \iota^*\h_S(X)
  \end{equation*}
is an $\A^1$-local equivalence of $\sC$-fibred spaces.
By universality of colimits it suffices to show that, for any $Y\in\sC_{/S}$ and any morphism $\varphi : \h_S(Y) \to \iota^*\h_S(X)$ (corresponding to a morphism $Y \to X$ in $\sD_{/S}$), the base change
  \begin{equation*}
    \iota^*\h_S(X \times \A^1) \fibprod_{\iota^*\h_S(X)} \h_S(Y) \to \h_S(Y)
  \end{equation*}
is an $\A^1$-local equivalence.
Since the morphism $\varphi$ factors as $\h_S(Y) \to \iota^*\h_S(Y) \to \iota^*\h_S(X)$, the morphism in question is a base change of the morphism
  \begin{equation*}
    \iota^*\h_S(X \times \A^1) \fibprod_{\iota^*\h_S(X)} \iota^*\h_S(Y) \to \iota^*\h_S(Y),
  \end{equation*}
which itself is identified with the canonical morphism
  \begin{equation*}
    \h_S(Y \times \A^1) \to \h_S(Y),
  \end{equation*}
since $\iota^*$ and $\h_S$ commute with limits and $\iota^*\iota_! = \id$.
This is an $\A^1$-local equivalence, so the claim follows.
\end{proof}

\begin{proof}[Proof of \propref{prop:MotSpc vs MOTSPC}]
Since $\iota_!$ preserves motivic equivalences (\lemref{lem:iota and local equivs}), its right adjoint $\iota^*$ preserves motivic spaces and induces a functor $\iota^* : \MotSpc(\sD_{/S}) \to \MotSpc(\sC_{/S})$, right adjoint to $\L\iota_!$.
Similarly, \lemref{lem:iota and local equivs} also implies that the right Kan extension functor $\iota_*$ preserves motivic spaces and defines a right adjoint to $\iota^* : \MotSpc(\sD_{/S}) \to \MotSpc(\sC_{/S})$.
Now the fully faithfulness of $\L\iota_!$, which is equivalent to invertibility of the unit map $\iota^*\L\iota_! \to \id$, follows by passage to left adjoints from the fully faithfulness of $\iota_*$ (which is equivalent to invertibility of the counit map $\iota^*\iota_*\to\id$).
The description of the essential image follows from \cite[Lem.~5.5.8.14]{HTT}.
\end{proof}

\begin{cor}\label{cor:Lhtp commutes with iota^*}
  There is a canonical invertible natural transformation
  \begin{equation*}
    \Lhtp \iota^* \to \iota^* \Lhtp.
  \end{equation*}
\end{cor}

\begin{proof}
It follows from \lemref{lem:iota and local equivs} that $\iota^* \Lhtp$ takes $\A^1$-invariant values, so the natural transformation in question is induced by the canonical map $\id \to \Lhtp$.
The fact that it is invertible follows from the formula \eqref{eq:formula for Lhtp}, which is valid for both $\sC$- and $\sD$-fibred spaces.
\end{proof}


\ssec{Functoriality}
\label{ssec:nil/functoriality}

We now record the various functorialities of $\sC$-fibred motivic spaces as the base varies; this works exactly as in the case $\sC = \Sm$ treated in \cite[Subsect.~2.5]{KhanLocalization}.
We then discuss the compatibility of these operations, as well as products and internal homs (\remref{rem:i^*Hom}), under the operation of extension along an inclusion of admissible subcategories (\propref{prop:MotSpc vs MOTSPC}).

\begin{notat}\label{notat:pasdfo}
  Let $f : T \to S$ be a morphism of affine spectral schemes.
  Let $\sC_{/S} \subseteq \Aff_{/S}$ be an admissible subcategory, and choose also an admissible subcategory $\sC_{/T} \subseteq \Aff_{/T}$ which contains the base changes $X \fibprod_S T$ of every $X \in \sC_{/S}$.
  A minimal such can be constructed as in \defnref{defn:broad}.
\end{notat}

\begin{constr}
  Under the notation of \ref{notat:pasdfo}, the base change functor $\Aff_{/S} \to \Aff_{/T}$ restricts to $\sC_{/S} \to \sC_{/T}$.
  \begin{defnlist}
    \item
    The direct image functor $f_*$ on $\sC$-fibred spaces is given by restriction along the base change functor $\sC_{/S} \to \sC_{/T}$.
    The latter preserves Nisnevich covering families and $\A^1$-projections, so $f_*$ preserves motivic spaces.
    Its left adjoint $f^*$ on motivic spaces is characterized uniquely by commutativity with colimits and the formula $f^*_{\sD}(\L\hspc[S]{X}) \simeq \L\hspc[T]{X \fibprod_S T}$ for $X \in \sC_{/S}$.

    \item\label{item:f_sharp^SPC}
    Suppose that the morphism $f : T \to S$ exhibits $T$ as an object of the full subcategory $\sC_{/S} \subseteq \Aff_{/S}$.
    Then the base change functor $\sC_{/S} \to \sC_{/T}$ admits a left adjoint, the forgetful functor
    \[ (X \to T) \mapsto (X \to T \xrightarrow{f} S) \]
    which preserves Nisnevich covering families and $\A^1$-projections.
    In this case $f^*$ is given by restriction along this forgetful functor, hence preserves motivic spaces and admits a left adjoint $f_\sharp$ characterized uniquely by commutativity with colimits and the formula $f_\sharp^{\sD}(\L\hspc[T]{X}) \simeq \L\hspc[S]{X}$ for $X \in\sC_{/T}$.
  \end{defnlist}
\end{constr}

We now discuss the compatibility of these operations under the embedding of \propref{prop:MotSpc vs MOTSPC}.
For this we fix the following notation.

\begin{notat}\label{notat:pjeqpkm}
  Let $f : T \to S$ be a morphism of affine spectral schemes.
  Fix an inclusion of admissible subcategories $\sC_{/S} \sub \sD_{/S}$ as in \notatref{notat:inclusion}.
  Fix similarly an inclusion $\sC_{/T} \sub \sD_{/T}$ of admissible subcategories both satisfying the condition of \notatref{notat:pasdfo}.
  Write $\MotSpc(\sC_{/S})$, $\MotSpc(\sD_{/S})$, $\MotSpc(\sC_{/T})$, and $\MotSpc(\sD_{/T})$ for the \inftyCats of motivic fibred spaces formed with respect to these choices.
\end{notat}

\begin{rem}\label{rem:functoriality and iota}
  The condition of \notatref{notat:pasdfo} guarantees that the base change functor $X \mapsto X \fibprod_S T$ commutes with the inclusions $\iota : \sC_{/S} \to \sD_{/S}$ and $\iota : \sC_{/T} \to \sD_{/T}$.
  From this it follows that the functor $f_*$ commutes with $\iota^*$ and that $f^*$ commutes with $\L\iota_!$.
  That is, we have commutative squares
  \[
    \begin{tikzcd}
      \MotSpc(\sC_{/S}) \ar{r}{\L\iota_!}\ar{d}{f^*}
        & \MotSpc(\sD_{/S}) \ar{d}{f^*}
      \\
      \MotSpc(\sC_{/T}) \ar{r}{\L\iota_!}
        & \MotSpc(\sD_{/T}),
    \end{tikzcd}
    \qquad
    \begin{tikzcd}
      \MotSpc(\sD_{/T}) \ar{r}{\iota^*}\ar{d}{f_*}
        & \MotSpc(\sC_{/T}) \ar{d}{f_*}
      \\
      \MotSpc(\sD_{/S}) \ar{r}{\iota^*}
        & \MotSpc(\sC_{/S}).
    \end{tikzcd}
  \]
  Similarly, if $f$ exhibits $T$ as an object of $\sC_{/S}$, then $f_\sharp$ commutes with $\L\iota_!$ and $f^*$ commutes with $\iota^*$.
\end{rem}

The following compatibility is less obvious:

\begin{prop}\label{prop:i_* and L iota_!}
  With notation as in \ref{notat:pjeqpkm}, assume that $\sC_{/S}$ is \emph{narrow}.
  Let $i : Z \hook S$ be a closed immersion of affine spectral schemes with affine open complement.
  Then there is a canonical invertible natural transformation
  \begin{equation}\label{eq:i_* and L iota_!}
    \L\iota_!\circ i_*^{\sC} \to i_*^{\sD}\circ\L\iota_!
  \end{equation}
  of functors $\MotSpc(\sC_{/Z}) \to \MotSpc(\sD_{/S})$, where the decorations indicate whether the functor is defined on $\sC$-fibred or $\sD$-fibred motivic spaces.
\end{prop}

The $\Sm$-fibred localization theorem \cite[Thm.~3.2.2]{KhanLocalization} implies the same for $\sC$-fibred motivic spaces (for any narrow $\sC$).
We will deduce \propref{prop:i_* and L iota_!} by combining this with the following $\sD$-fibred variant:

\begin{thm}[Localization] \label{thm:localization}
  Let the notation be as in \propref{prop:i_* and L iota_!}, and let $j : U \to S$ be the open immersion complementary to $i$.
  Let $\sF \in \MotSpc(\sD_{/S})$ be a $\sD$-fibred motivic space over $S$.
  If $\sF$ belongs to the essential image of the functor $\L\iota_! : \MotSpc(\sC_{/S}) \to \MotSpc(\sD_{/S})$ (\propref{prop:MotSpc vs MOTSPC}), then there is a cocartesian square
  \begin{equation*}
    \begin{tikzcd}
      j_\sharp j^* (\sF) \arrow{r}\arrow{d} & \sF \arrow{d} \\
      j_\sharp (\pt_U) \arrow{r} & i_* i^* (\sF)
    \end{tikzcd}
  \end{equation*}
  of $\sD$-fibred motivic spaces over $S$.
\end{thm}

\begin{proof}
From \cite[Prop.~3.1.4]{KhanLocalization} it follows that the functor $i_* : \MotSpc(\sD_{/Z}) \to \MotSpc(\sD_{/S})$ commutes with contractible\footnotemark~colimits (exactly as in the proof of Thm.~3.1.1 in \emph{loc. cit}).
By \propref{prop:MotSpc vs MOTSPC} we may assume that $\sF$ is the motivic localization of some $X \in \sC_{/S}$.
Since $\sC_{/S}$ is narrow, we can moreover assume by \examref{exam:narrow=Sm-fibred} that it belongs to $\sA^0_{/S}$ (\examref{exam:minimal adm}), i.e., that it admits an étale morphism to some $\A^n_S$, $n\ge 0$.
Then we can proceed exactly as in \cite[Subsect.~4.3]{KhanLocalization}.
\end{proof}

\begin{proof}[Proof of \propref{prop:i_* and L iota_!}]
The natural transformation \eqref{eq:i_* and L iota_!} is the composite
  \begin{equation*}
    \L\iota_!\circ i_*^{\sC}
      \xrightarrow{\mrm{unit}} i_*^{\sD} i^*_{\sD} \L\iota_!\circ i_*^{\sC}
      \simeq i_*^{\sD} \L\iota_!\circ i^*_{\sC} i_*^{\sC}
      \xrightarrow{\mrm{counit}} i_*^{\sD}\circ\L\iota_!,
  \end{equation*}
where the identification $\L \iota_! \circ i^*_{\sC} \simeq i^*_{\sD} \circ \L \iota_!$ comes from \remref{rem:functoriality and iota}.
By \cite[Prop.~2.4.4]{KhanLocalization} it will suffice to show that the canonical morphism
  \begin{equation*}
    \L \iota_!\circ i_*^{\sC}(\L \h_Z(X)) \to i_*^{\sD} (\L \h_Z(X)) \\
  \end{equation*}
is invertible for all $X \in \sC_{/Z}$.
Using \cite[Prop.~3.1.4]{KhanLocalization}, we may assume that $X$ is of the form $Y \fibprod_S Z$ for some $Y\in\sC_{/S}$.
Then we conclude by comparing the description of $i_*^{\sC}i^*_{\sC}(\L \h_S(Y))$ given by \cite[Thm.~3.2.2]{KhanLocalization}, and the description of $i_*^{\sD}i^*_{\sD}(\L \h_S(Y))$ provided by \thmref{thm:localization}.
\end{proof}

Finally, we discuss the compatibility of the functors $\L\iota_!$ and $\iota^*$ with products and internal homs.
Note that, just as in the $\sC$-fibred case \cite[Rem.~2.4.2]{KhanLocalization}, the full subcategory $\MotSpc(\sD_{/S}) \subseteq \Spc(\sD_{/S})$ is closed under formation of internal homs.

\begin{rem}\label{rem:i^*Hom}\leavevmode
  \begin{defnlist}
    \item
    Since the functor $\iota^* : \MotSpc(\sD_{/S}) \to \MotSpc(\sC_{/S})$ preserves limits (see proof of \propref{prop:MotSpc vs MOTSPC}), it is symmetric monoidal with respect to the cartesian product.
    
    \item\label{item:Li! monoidal}
    The functor $\L\iota_! : \MotSpc(\sC_{/S}) \to \MotSpc(\sD_{/S})$ is also symmetric monoidal.
    Indeed, since $\L$ preserves finite products by \remref{rem:SPC localization basics}(iii), it suffices to show the claim for $\iota_! : \Spc(\sC_{/S}) \to \Spc(\sD_{/S})$.
    For this we may reduce to representables which is obvious.

    \item
    For any $\sC$-fibred motivic space $\sF \in \MotSpc(\sD_{/S})$ and $\sD$-fibred motivic space $\sG \in \MotSpc(\sD_{/S})$, there is a canonical isomorphism
    \[
      \iota^* \uHom(\L\iota_!(\sF), \sG) \to \uHom(\sF, \iota^*(\sG))
    \]
    of $\sC$-fibred motivic spaces, where $\uHom$ is taken in $\MotSpc(\sD_{/S})$ on the left and in $\MotSpc(\sC_{/S})$ on the right.
    This follows by adjunction from \itemref{item:Li! monoidal}.
  \end{defnlist}
\end{rem}


\ssec{Classical fibred spaces}
\label{ssec:nil/classical}

In this subsection we set up, for a classical affine scheme $S$, a classical variant of the \inftyCat of $\sC$-fibred motivic spaces.
For $S$ a spectral affine scheme we then define a pair of adjoint functors
\[
  \L v_! : \MotSpc(\sC_{/S}) \to \MotSpc(\sC^\cl_{/S_\cl}),\quad
  v^* : \MotSpc(\sC^\cl_{/S_\cl}) \to \MotSpc(\sC_{/S}).
\]
Later we will focus on understanding these adjunctions when $\sC_{/S}$ is broad (\ssecref{ssec:nil/nil-localization}) and narrow (\ssecref{ssec:nil/comparison}).

The following is a classical analogue of Definitions~\ref{defn:adm} and \ref{defn:narrow}.

\begin{defn}\label{defn:adm cl}
  Let $S$ be an affine scheme.
  Denote by $\AffCl_{/S}$ the category of classical affine schemes over $S$.
  We say that a full subcategory $\sC_{/S} \sub \AffCl_{/S}$ is \emph{admissible} if it is essentially small and satisfies the following conditions:
  \begin{defnlist}[series=adm]
    \item
    The affine scheme $S$ (viewed over $S$ via the identity) belongs to $\sC_{/S}$.
    \item
    If $X$ belongs to $\sC_{/S}$ and $Y$ is étale over $X$, then $Y$ belongs to $\sC_{/S}$.
    \item
    If $X$ belongs to $\sC_{/S}$, then $X \times \A^n_\cl$ belongs to $\sC_{/S}$ for every $n\ge 0$.
  \end{defnlist}

  For example, the full subcategory $\SmCl_{/S} \sub \AffCl_{/S}$ of \emph{smooth} affine schemes over $S$ (where smoothness is understood in the sense of classical algebraic geometry) is admissible.
  We say that an admissible subcategory $\sC_{/S} \sub \AffCl_{/S}$ is \emph{narrow} if it is contained in $\SmCl_{/S}$.
\end{defn}

\begin{exam}\label{exam:A^cl}
  Let $S$ be a \emph{spectral} affine scheme and $\sC_{/S} \sub \Aff_{/S}$ an admissible subcategory.
  Then the full subcategory $\sC^\cl_{/S_\cl}\sub\AffCl_{/S_\cl}$ spanned by the classical truncations $X_\cl$ of all objects $X \in \sC_{/S}$ is admissible.
  The first condition follows from \defnref{defn:adm}\ref{item:adm1}, the second follows from \defnref{defn:adm}\ref{item:adm2} and \cite[Thm.~7.5.0.6]{HA-20170918}, and the third follows from \defnref{defn:adm}\ref{item:adm3} and the fact that $\A^n_\cl \simeq (\A^n)_\cl$.
\end{exam}

\begin{defn}\label{defn:SmCl-fibred space}
  Let $S$ be an affine scheme and let $\sC^\cl_{/S} \sub \AffCl_{/S}$ be an admissible subcategory.
  A \emph{$\sC^\cl$-fibred space} over $S$ is a presheaf of spaces on $\sC^\cl_{/S}$.
  We say that a $\sC^\cl$-fibred space $\sF$ over $S$ satisfies \emph{Nisnevich excision} if it is reduced, and for any $X \in \sC^\cl_{/S}$ and any Nisnevich square $Q$ over $X$, the induced square of spaces $\Gamma(Q, \sF)$ is cartesian.
  We say that $\sF$ satisfies \emph{$\A^1_\cl$-homotopy invariance} if for any $X \in \sC^\cl_{/S}$, the canonical map $\Gamma(X, \sF) \to \Gamma(X\times\A^1_\cl,\sF)$ is invertible, where $\A^1_\cl = \Spec(\bZ[T]) \simeq (\A^1)_\cl$ denotes the classical affine line.
  A \emph{motivic $\sC^\cl$-fibred space} is a $\sC^\cl$-fibred space that satisfies Nisnevich excision and $\A^1_\cl$-homotopy invariance.
\end{defn}

\begin{exam}
  Let $S = \Spec(R)$ be a spectral affine scheme.
  Then the \inftyCat of motivic $\SmCl$-fibred spaces $\MotSpc(\SmCl_{/S_\cl})$ is equivalent to $\MotSpc^\cl(\pi_0(R))$ as defined in \thmref{thm:intro/comp1}.
\end{exam}

For the remainder of this subsection, we fix the following notation:

\begin{notat}
  Let $S$ be a spectral affine scheme.
  Fix an admissible subcategory $\sC_{/S} \sub \Aff_{/S}$ (\defnref{defn:adm}).
  Let $\sC^\cl_{/S}$ be the induced admissible subcategory of $\AffCl_{/S}$ as in \examref{exam:A^cl}.
  We denote by
  \[ \Spc(\sC_{/S}), \quad \text{resp.}~\Spc(\sC^\cl_{/S_\cl}), \]
  the \inftyCat of $\sC$-fibred spaces over $S$, resp. of $\sC^\cl$-fibred spaces over $S_\cl$, and by
  \[ \MotSpc(\sC_{/S}), \quad \text{resp.}~\MotSpc(\sC^\cl_{/S_\cl}), \]
  the full subcategory of motivic objects.
\end{notat}

\begin{constr}\label{constr:v}
  The operation of passing to classical truncations,
  \[ (X \to S) \mapsto (X_\cl \to S_\cl), \]
  defines a canonical functor $v : \Aff_{/S} \to \AffCl_{/S_\cl}$ which restricts to $v : \sC_{/S} \to \sC^\cl_{/S_\cl}$.
  Denote by $v^* : \Spc(\sC^\cl_{/S_\cl}) \to \Spc(\sC_{/S})$ the functor of restriction along $v$, and by $v_! : \Spc(\sC_{/S}) \to \Spc(\sC^\cl_{/S_\cl})$ its left adjoint given by left Kan extension of $v$.
  Recall that $v_!$ is uniquely characterized by commutativity with colimits and the formula $v_!(\hspc[S]{X}) \simeq \hspc[S_\cl]{X_\cl}$ for all $X \in \sC_{/S}$.
\end{constr}

\begin{lem}\label{lem:v_!,v^* and motivic equivalences}\leavevmode
  \begin{thmlist}
    \item
    The functor $v_!$ preserves Nisnevich-local equivalences, $\A^1$-local equivalences, and motivic equivalences.

    \item
    The functor $v^*$ preserves Nisnevich excisive spaces and Nisnevich-local equivalences.
    In particular, it commutes with $\LNis$; that is, there is a canonical invertible natural transformation $\LNis v^* \to v^* \LNis$.

    \item
    The functor $v^*$ sends $\A^1_\cl$-invariant $\sC^\cl$-fibred spaces to $\A^1$-invariant $\sC$-fibred spaces.
    In particular, it sends motivic $\sC^\cl$-fibred spaces to motivic $\sC$-fibred spaces.
  \end{thmlist}
\end{lem}

In particular, we find that the functors $v_!$ and $v^*$ descend to a pair of adjoint functors
\begin{equation}\label{eq:Lv_!}
  \L v_! : \MotSpc(\sC_{/S}) \to \MotSpc(\sC^\cl_{/S_\cl}),
  \quad v^* : \MotSpc(\sC^\cl_{/S_\cl}) \to \MotSpc(\sC_{/S}).
\end{equation}

\begin{proof}
  The first claim follows from the fact that $v$ preserves Nisnevich squares and sends $\A^1$ to $\A^1_\cl$.
  By adjunction it follows that $v^*$ preserves Nisnevich-excisive spaces and sends $\A^1_\cl$-invariant spaces to $\A^1$-invariant spaces.
  It remains to show that $v^*$ preserves Nisnevich-local equivalences.
  For this it is sufficient to check that the functor $v$ is cocontinuous for the Nisnevich topology, i.e., that for all $X \in \sC_{/S}$, any Nisnevich covering of $X_\cl$ lifts to a Nisnevich covering of $X$ (cf.~\cite[Def.~3.1.5]{KhanLocalization}).
  This follows from \cite[Thm.~7.5.0.6]{HA-20170918}.
\end{proof}


\ssec{Nil-localization}
\label{ssec:nil/nil-localization}

In this subsection we study the adjunction \eqref{eq:Lv_!} when the admissible subcategory $\sC_{/S} \sub \Aff_{/S}$ is \emph{broad}.
In this case, we find that the classical construction $\MotSpc(\sC^\cl_{/S_\cl})$ is a left Bousfield localization of the spectral variant $\MotSpc(\sC_{/S})$ (\thmref{thm:nil-localization}).

\begin{notat}
  Let $S$ be a spectral affine scheme.
  Fix a broad subcategory $\sB_{/S} \sub \Aff_{/S}$ (\defnref{defn:broad}), and let $\sB^\cl_{/S}$ be the induced admissible subcategory of $\AffCl_{/S}$ as in \examref{exam:A^cl}.
  Consider the \inftyCats
  \[ \MotSpc(\sB_{/S}) \sub \Spc(\sB_{/S}), \quad \MotSpc(\sB^\cl_{/S_\cl}) \sub \Spc(\sB^\cl_{/S_\cl}). \]
\end{notat}

\begin{defn}\label{defn:nil-local}
  A $\sB$-fibred space $\sF \in \Spc(\sB_{/S})$ is called \emph{nil-local} if for any $X \in \sB_{/S}$, the canonical map of spaces
  \begin{equation*}
    \Gamma(X, \sF) \to \Gamma(X_\cl, \sF)
  \end{equation*}
  is invertible.
\end{defn}

\begin{thm}\label{thm:nil-localization}
The functor $v^* : \MotSpc(\sB^\cl_{/S_\cl}) \to \MotSpc(\sB_{/S})$ is fully faithful, and induces an equivalence
  \begin{equation*}
    v^* : \MotSpc(\sB^\cl_{/S_\cl}) \to \MotSpc_\nil(\sB_{/S})
  \end{equation*}
from the \inftyCat of $\sB^\cl$-fibred motivic spaces over $S_\cl$ to the \inftyCat $\MotSpc_\nil(\sB_{/S}) \sub \MotSpc(\sB_{/S})$ of nil-local $\sB$-fibred motivic spaces over $S$.
In particular, the functor $\L v_! : \MotSpc(\sB_{/S}) \to \MotSpc(\sB^\cl_{/S_\cl})$ is a left Bousfield localization.
\end{thm}

The proof of \thmref{thm:nil-localization} relies on an analysis of the behaviour of the functors
\[
  v_! : \Spc(\sB_{/S}) \to \Spc(\sB^\cl_{/S_\cl}),
  \quad
  v^* : \Spc(\sB^\cl_{/S_\cl}) \to \Spc(\sB_{/S}) \]
with respect to $\A^1$-local and Nisnevich-local equivalences, specializing \lemref{lem:v_!,v^* and motivic equivalences} to the broad case:

\begin{prop}\label{prop:v_!,v^* and motivic equivalences}\leavevmode
  \begin{thmlist}
    \item\label{item:v_!Nis}
    The functor $v_!$ preserves Nisnevich excisive spaces and Nisnevich-local equivalences.
    In particular, it commutes with $\LNis$; that is, there is a canonical invertible natural transformation $\LNis v_! \to v_! \LNis$.

    \item\label{item:v^*Nis}
    The functor $v^*$ preserves Nisnevich excisive spaces and Nisnevich-local equivalences.
    In particular, it commutes with $\LNis$; that is, there is a canonical invertible natural transformation $\LNis v^* \to v^* \LNis$.

    \item\label{item:v_!A1}
    The functor $v_!$ sends $\A^1$-local equivalences to $\A^1_\cl$-local equivalences.

    \item\label{item:v^*A1}
    The canonical natural transformations
      \begin{equation*}
        \Lhtp v^* \to v^* \Lhtpcl
        \quad\textit{and}\quad
        v_! \Lhtp v^* \to \Lhtpcl
      \end{equation*}
    are invertible.

    \item\label{item:v_!mot}
    The functor $v_!$ preserves motivic equivalences.

    \item\label{item:v^*mot}
    The canonical natural transformations
      \begin{equation*}
        \L v^* \to v^* \L
        \quad\textit{and}\quad
        v_! \L v^* \to \L
      \end{equation*}
    are invertible.
  \end{thmlist}
\end{prop}

The key feature of the broad case is the existence of a left adjoint $u$ to the functor $v : \sB_{/S} \to \sB^\cl_{/S_\cl}$:

\begin{rem}\label{rem:u}
  Any affine scheme $X_0$ over $S_\cl$ can be viewed as a discrete affine spectral scheme over $S$ (by composition with the canonical morphism $S_\cl \to S$).
  This defines a functor $u : \AffCl_{/S_\cl} \to \Aff_{/S}$, left adjoint to $v$.
  Note that $u$ is fully faithful (as the unit map $X_0 \to vu(X_0)$ is always invertible).
  Since $\sB$ is \emph{broad}, our choice of $\sB^\cl_{/S_\cl}$ (\examref{exam:A^cl}) guarantees that $u$ restricts to a functor $u : \sB^\cl_{/S_\cl} \to \sB_{/S}$.
  Restriction along $u$ defines a functor $u^* : \Spc(\sB_{/S}) \to \Spc(\sB^\cl_{/S_\cl})$, which admits fully faithful left and right adjoints $u_!$ and $u_*$, respectively.
  By adjunction, we have identifications $v_! \simeq u^*$ and $v^* \simeq u_*$.
  In particular, it follows that the co-unit
  \[ v_!v^* \to \id \]
  is invertible.
\end{rem}

\begin{rem}\label{rem:v^*(F) is nil-local}
  Note that a $\sB$-fibred space $\sF \in \Spc(\sB_{/S})$ is nil-local if and only if it belongs to the essential image of $v^*$, or equivalently if and only if the unit map $\sF \to v^*v_!(\sF) \simeq v^*u^*(\sF)$ is invertible.
  Indeed the counit map $uv(X) \to X$ is the inclusion of the classical truncation, so the map $\Gamma(X, \sF) \to \Gamma(X_\cl, \sF)$ is canonically identified with the map
  \[ \Gamma(X, \sF) \to \Gamma(X, v^*u^*(\sF)) \]
  for every $X \in \sB_{/S}$.
\end{rem}

\begin{proof}[Proof of \propref{prop:v_!,v^* and motivic equivalences}]
  We already know from \lemref{lem:v_!,v^* and motivic equivalences} that $v_!$ sends Nisnevich-local equivalences to Nisnevich-local equivalences, $\A^1$-local equivalences to $\A^1_\cl$-local equivalences, and motivic equivalences to motivic equivalences.
  We also know that its right adjoint $v^*$ sends Nisnevich excisive spaces to Nisnevich excisive spaces, $\A^1_\cl$-invariant spaces to $\A^1$-invariant spaces, and motivic spaces to motivic spaces.

  Let $u : \sB^\cl_{/S_\cl} \to \sB_{/S}$ be as in \remref{rem:u}, so that $v_! \simeq u^*$.
  Since $u$ preserves Nisnevich squares, the functor $u_!$ preserves Nisnevich-local equivalences.
  Hence its right adjoint $u^* \simeq v_!$ preserves Nisnevich excisive spaces.
  This proves claim~\ref{item:v_!Nis}.

  Consider claim~\ref{item:v^*A1}.
  Since $v^*$ preserves $\A^1$-invariant spaces, the natural transformation $\id \to \Lhtpcl$ induces a transformation
  \begin{equation}\label{eq:prop:v_!,v^* and motivic equivalences (iv)}
    \Lhtp v^* v_! \to \Lhtp v^* \Lhtpcl v_! \simeq v^* \Lhtpcl v_!
  \end{equation}
  which we claim is invertible.
  For every $X \in \sB_{/S}$, let $\sI_X \sub \Aff_{/X}$ denote the non-full subcategory whose objects are spectral affine spaces $\A^n \times X$ over $X$, and whose morphisms are projections.
  A variant of the formula \eqref{eq:formula for Lhtp} (see \cite[Prop.~3.4]{HoyoisEquivariant}) then yields the functorial isomorphisms
    \begin{align*}
      \Gamma\big(X, \Lhtp v^*u^*(\sF)\big)
        &\simeq \colim_{\sI_X^\op} \Gamma\big(\A^n \times X, v^*u^*(\sF)\big)\\
        &\simeq \colim_{\sI_X^\op} \Gamma\big(\A^n_\cl \times X_\cl, \sF\big),
    \end{align*}
  by \remref{rem:v^*(F) is nil-local}.
  Similarly, if we write $\sI^\cl_{X_\cl} \sub \AffCl_{/X_\cl}$ for the subcategory of classical affine spaces $\A^n_\cl \times X_\cl$ and projections between them, then \cite[Prop.~3.4]{HoyoisEquivariant} again yields
    \begin{align*}
      \Gamma\big(X, v^*\Lhtpcl u^*(\sF)\big)
        &\simeq \Gamma\big(X_\cl, \Lhtpcl u^*(\sF)\big)\\
        &\simeq \colim_{(\sI^\cl_{X_\cl})^\op} \Gamma\big(\A^n_\cl \times X_\cl, u^*(\sF)\big)\\
        &\simeq \colim_{(\sI^\cl_{X_\cl})^\op} \Gamma\big(\A^n_\cl \times X_\cl, \sF\big).
    \end{align*}
  Since $v : \sB_{/X} \to \sB^\cl_{/X_\cl}$ induces an equivalence $\sI_X \simeq \sI^\cl_{X_\cl}$, it follows that \eqref{eq:prop:v_!,v^* and motivic equivalences (iv)} is invertible.
  Applying $v^*$ on the right, we deduce that the canonical transformation $\Lhtp v^* \to v^* \Lhtpcl$ is also invertible (since $v_!v^* \simeq \id$ by \remref{rem:u}).
  Applying $v_!$ on the left, we also obtain the invertible transformation $v_! \Lhtp v^* \to \Lhtpcl$.

  Claim~\ref{item:v^*mot} follows from claims~\ref{item:v^*Nis} and \ref{item:v^*A1} in view of the formula \eqref{eq:description of Lmot} (and the analogous formula for $\sB^\cl$-fibred spaces).
\end{proof}

\begin{proof}[Proof of \thmref{thm:nil-localization}]
  By \remref{rem:u} we know that the functor \[v^* : \MotSpc(\sB^\cl_{/S_\cl}) \to \MotSpc(\sB_{/S})\] is fully faithful.
  Its essential image is spanned by objects $\sF \in \MotSpc(\sB_{/S})$ for which the unit map $\sF \to v^*\L v_!(\sF)$ is invertible.
  By \remref{rem:v^*(F) is nil-local}, this condition implies that $\sF$ is nil-local.
  Conversely if $\sF$ is nil-local, so that the unit map $\sF \to v^*v_!(\sF)$ is invertible (again by \remref{rem:v^*(F) is nil-local}), then using \propref{prop:v_!,v^* and motivic equivalences}\ref{item:v^*mot} we see that the induced map $\sF \simeq \L(\sF) \to \L v^*v_!(\sF) \simeq v^*\L v_!(\sF)$ is also invertible.
\end{proof}


\ssec{Nil descent}
\label{ssec:nil/nil-local}

\begin{notat}\label{notat:narrow broad}
  Let $S$ be a spectral affine scheme.
  Fix a narrow subcategory $\sA_{/S} \sub \Aff_{/S}$ and a broad subcategory $\sB_{/S} \sub \Aff_{/S}$ containing $\sA_{/S}$.
  We consider the \inftyCats
  \[ \MotSpc(\sA_{/S}) \sub \Spc(\sA_{/S}), \quad \MotSpc(\sB_{/S}) \sub \Spc(\sB_{/S}).\]
  Let $\iota : \sA_{/S} \to \sB_{/S}$ denote the inclusion, so that we have the fully faithful functor (see \propref{prop:MotSpc vs MOTSPC})
  \[
    \L\iota_! : \MotSpc(\sA_{/S}) \to \MotSpc(\sB_{/S}).
  \]
\end{notat}

We are now in a position to state and prove the following result:

\begin{thm}\label{thm:nil descent}
  Let $\sF \in \MotSpc(\sA_{/S})$ be an $\sA$-fibred motivic space.
  Then the $\sB$-fibred motivic space $\L\iota_!(\sF) \in \MotSpc(\sB_{/S})$ is nil-local (\defnref{defn:nil-local}).
\end{thm}

\begin{proof}
  Set $\sF^+ = \L\iota_!(\sF)$.
  Let $X\in\sB_{/S}$ with structural morphism $f : X \to S$ and choose subcategories $\sA_{/X} \sub \sB_{/X}$ of $\Aff_{/X}$, narrow and broad, respectively, and both satisfying the condition of \notatref{notat:pasdfo}.
  By adjunction, there are canonical isomorphisms
    \begin{align*}
      \Gamma(X, \sF^+)
        &\simeq \Maps(\pt_X, f^*_{\sB}(\sF^+)),\\
      \Gamma(X_\cl, \sF^+)
        &\simeq \Maps(\h_X(X_\cl), f^*_{\sB}(\sF^+))\\
        &\simeq \Maps(i_\sharp^\sB i^*_\sB (\pt_X), f^*_{\sB}(\sF^+))\\
        &\simeq \Maps(\pt_X, i_*^\sB i^*_\sB f^*_{\sB}(\sF^+)),
    \end{align*}
  where all the mapping spaces are formed in $\MotSpc(\sB_{/X})$.
  Under these identifications the map $\Gamma(X, \sF^+) \to \Gamma(X_\cl, \sF^+)$ is induced by the unit morphism
    \begin{equation*}
      f^*_{\sB}(\sF^+) \to i_*^{\sB} i^*_{\sB} f^*_{\sB}(\sF^+)
    \end{equation*}
  in $\MotSpc(\sB_{/X})$.
  By \remref{rem:functoriality and iota} and \propref{prop:i_* and L iota_!} this morphism is the image by $\L\iota_!$ of the unit morphism
    \begin{equation*}
      f^*_{\sA} (\sF) \to i_*^{\sA} i^*_{\sA} f^*_{\sA} (\sF)
    \end{equation*}
  in $\MotSpc(\sA_{/X})$.
  Since $i$ is a closed immersion with empty complement, this morphism is invertible by the nilpotent invariance property of $\MotSpc(\sA_{/X})$, see \cite[Cor.~3.2.7]{KhanLocalization} (which applies to any narrow subcategory and not just $\Sm_{/X}$).
\end{proof}


\ssec{The comparison}
\label{ssec:nil/comparison}

\begin{notat}\label{notat:narrow broad classical}
  Let $S$ be a spectral affine scheme.
  We again fix narrow and broad subcategories $\sA_{/S} \sub \Aff_{/S}$ and $\sB_{/S} \sub \Aff_{/S_\cl}$ as in \notatref{notat:narrow broad}.
  We let $\sA^\cl_{/S}$ and $\sB^\cl_{/S}$ be the induced admissible subcategories of $\AffCl_{/S}$ as in \examref{exam:A^cl}.
  To simplify notation set
  \[ \SPC(S) := \Spc(\sB_{/S}), \quad \SPC^\cl(S_\cl) := \Spc(\sB^\cl_{/S_\cl}), \quad \MOTSPC(S) := \MotSpc(\sB_{/S}), \quad \MOTSPC^\cl(S_\cl) := \MotSpc(\sB^\cl_{/S_\cl}), \]
  and similarly
  \[ \Spc(S) := \Spc(\sA_{/S}), \quad \Spc^\cl(S_\cl) := \Spc(\sA^\cl_{/S_\cl}), \quad \MotSpc(S) := \MotSpc(\sA_{/S}), \quad \MotSpc^\cl(S_\cl) := \MotSpc(\sA^\cl_{/S_\cl}). \]
  Recall from \examref{exam:narrow=Sm-fibred} that this notation agrees with that of \examref{exam:MotSpc(S)}, even though $\sA_{/S}$ is allowed to be any narrow subcategory.
\end{notat}

Let $v : \sB_{/S} \to \sB_{/S_\cl}$ and $w : \sA_{/S} \to \sA_{/S_\cl}$ be the classical truncation functors as in \constrref{constr:v}.
In this subsection we will prove the following result, which gives the equivalence between (i) and (iii) in \thmref{thm:intro/comp1}.

\begin{thm}\label{thm:comp}
  The adjunction of \eqref{eq:Lv_!},
  \[
    \L w_! : \MotSpc(S) \to \MotSpc^\cl(S_\cl),
    \quad
    w^* : \MotSpc^\cl(S_\cl) \to \MotSpc(S),
  \]
  is an equivalence of \inftyCats.
\end{thm}

The proof will combine the $\sB$-fibred nil-localization statement (\thmref{thm:nil-localization}) and nil descent for $\sA$-fibred spaces (\thmref{thm:nil descent}), as well as \propref{prop:MotSpc vs MOTSPC} and the following classical analogue of the latter:

\begin{prop}\label{prop:MotSpc^cl vs MOTSPC^cl}
  Consider the inclusion functor $\iota : \sA^\cl_{/S_\cl} \hook \sB^\cl_{/S_\cl}$.
  Let $\iota_! : \Spc^\cl(S_\cl) \to \SPC^\cl(S_\cl)$ denote the left Kan extension of $\iota$, left adjoint to the restriction functor $\iota^* : \SPC^\cl(S) \to \Spc^\cl(S)$.
  Then the assignment $\sF \mapsto \L\iota_!(\sF)$ induces a fully faithful functor of \inftyCats
    \begin{equation*}
      \L\iota_! : \MotSpc^\cl(S_\cl) \to \MOTSPC^\cl(S_\cl),
    \end{equation*}
  whose essential image is generated under sifted colimits by objects of the form $\L \h_{S_\cl}(X)$, where $X \in \sA^\cl_{/S_\cl}$ admits an étale $S_\cl$-morphism to a classical affine space $S_\cl \times \A^n_{\cl}$, for some $n\ge 0$.
\end{prop}
\begin{proof}
  Same proof as \propref{prop:MotSpc vs MOTSPC}.
\end{proof}

\begin{rem}\label{rem:w iota}
  Since $w : \sA_{/S} \to \sA^\cl_{/S_\cl}$ is the restriction of $v : \sB_{/S} \to \sB^\cl_{/S_\cl}$, we have commutative squares
  \begin{equation}\label{eq:v w squares}
    \begin{tikzcd}
      \MotSpc(S) \ar{r}{\L w_!}\ar{d}{\L\iota_!}
      & \MotSpc^\cl(S_\cl) \ar{d}{\L\iota_!}
      \\
      \MOTSPC(S) \ar{r}{\L v_!}
      & \MOTSPC^\cl(S_\cl),
    \end{tikzcd}
    \quad
    \begin{tikzcd}
      \MOTSPC^\cl(S_\cl) \ar{r}{v^*}\ar{d}{\iota^*}
      & \MOTSPC(S) \ar{d}{\iota^*}
      \\
      \MotSpc^\cl(S_\cl) \ar{r}{w^*}
      & \MotSpc(S).
    \end{tikzcd}
  \end{equation}
\end{rem}

\begin{proof}[Proof of \thmref{thm:comp}]
We show that the adjunction $(\L w_!, w^*)$ is an equivalence.
For any $\sA$-fibred motivic space $\sF \in \MotSpc(S)$, the $\sB$-fibred motivic space $\L \iota_!(\sF)$ is nil-local by \thmref{thm:nil descent}.
Therefore by \thmref{thm:nil-localization} the canonical map
\[ \L \iota_!(\sF) \to v^*\L v_! \L \iota_!(\sF) \]
is invertible.
Applying $\iota^*$ and using \propref{prop:MotSpc vs MOTSPC}, we deduce that the canonical map
\[ \sF \to \iota^* v^*\L v_! \L \iota_!(\sF) \]
is invertible.
This map is identified with the unit $\sF \to w^* \L w_!(\sF)$ under the identifications (\remref{rem:w iota} and \propref{prop:MotSpc^cl vs MOTSPC^cl})
\[
  \iota^* v^* \L v_! \L \iota_!
    \simeq w^* \iota^* \L \iota_! \L w_!
    \simeq w^* \L w_!.
\]
Since $\sF \in \MotSpc(S)$ was arbitrary, this shows that the unit
\[ \id \to w^* \L w_! \]
is invertible, hence $\L w_!$ is fully faithful.
It remains to show that $\MotSpc^\cl(S_\cl)$ is generated under colimits by objects of the form $\L \hspc[S_\cl]{X_\cl}$, where $X \in \sA_{/S}$.
But this follows from the definition of $\sA^\cl_{/S_\cl}$ (\examref{exam:A^cl}).
\end{proof}

The next few results are corollaries of Theorems~\ref{thm:nil-localization} and \ref{thm:comp}.

\begin{cor}\label{cor:comp asp}
  For any quasi-compact quasi-separated spectral algebraic space $S$, there are canonical equivalences of \inftyCats
  \[ \MotSpc(S) \simeq \MotSpc^\cl(S_\cl), \]
  where $\MotSpc(S)$ is as in \cite[Def.~2.4.1]{KhanLocalization} and $\MotSpc^\cl(S_\cl)$ its classical variant.
\end{cor}
\begin{proof}
  Classical truncation defines a functor $w$ from the \inftyCat of smooth spectral algebraic spaces over $S$ to the category of smooth classical algebraic spaces over $S_\cl$.
  This induces an adjunction
  \[
    \L w_! : \MotSpc(S) \to \MotSpc^\cl(S_\cl),
    \quad
    w^* : \MotSpc^\cl(S_\cl) \to \MotSpc(S)
  \]
  which globalizes that of \thmref{thm:comp}.
  To show that the unit and counit maps are invertible, we may use Nisnevich descent \cite[Prop.~2.5.7]{KhanLocalization} to reduce to the affine case proven in \thmref{thm:comp}.
\end{proof}

\begin{cor}\label{cor:Lmot(F) in terms of Lmot v_!(F)}
Let $S$ be an affine spectral scheme.
Then for any $\sA$-fibred space $\sF \in \Spc(S)$, the canonical map in $\MotSpc(S)$
  \begin{equation*}
    \L(\sF) \to w^*\L w_!(\sF)
  \end{equation*}
is invertible.
\end{cor}

\begin{proof}
By \thmref{thm:intro/comp1} the canonical map $\L(\sF) \isoto w^* \L w_! (\L(\sF))$ is invertible.
By \propref{prop:v_!,v^* and motivic equivalences} (v) the canonical map $w_!(\sF) \to w_! (\L(\sF))$ is a motivic equivalence, whence the claim.
\end{proof}

\begin{cor}\label{cor:iota_! commutes with v^*}
  Let $S$ be an affine spectral scheme.
  Then there is a commutative square
  \[
    \begin{tikzcd}
      \MotSpc^\cl(S_\cl) \ar{r}{w^*}\ar{d}{\L\iota_!}
      & \MotSpc(S) \ar{d}{\L\iota_!}
      \\
      \MOTSPC^\cl(S_\cl) \ar{r}{v^*}
      & \MOTSPC(S).
    \end{tikzcd}
  \]
\end{cor}

\begin{proof}
  The square is obtained by horizontally passing to right adjoints in the left-hand square in \eqref{eq:v w squares}, and thus commutes up to a natural transformation which we claim is invertible.
  Note that both clockwise and counterclockwise composites factor through the full subcategory $\MOTSPC_\nil(S)$ of nil-local objects by \thmref{thm:nil descent} and \remref{rem:v^*(F) is nil-local}.
  By \thmref{thm:nil-localization} it will therefore suffice to show that the natural transformation becomes invertible after post-composition with $\L v_! : \MOTSPC(S) \to \MOTSPC^\cl(S_\cl)$.
  This is immediate from the fact that $v^*$ is fully faithful (\thmref{thm:nil-localization}) and the commutativity of the right-hand square in \eqref{eq:v w squares}.
\end{proof}

\begin{cor}
  Let $S$ be an affine spectral scheme.
  Then there is a commutative square
  \[
    \begin{tikzcd}
      \MOTSPC_\nil(S) \ar{r}{\L v_!}\ar{d}{\iota^*}
      & \MOTSPC^\cl(S_\cl) \ar{d}{\iota^*}
      \\
      \MotSpc(S) \ar{r}{\L w_!}
      & \MotSpc^\cl(S_\cl).
    \end{tikzcd}
  \]
\end{cor}
\begin{proof}
  The square is obtained by horizontally passing to left adjoints in the right-hand square in \eqref{eq:v w squares} (and restricting to the full subcategory $\MOTSPC_\nil(S) \sub \MOTSPC(S)$), and thus commutes up to a natural transformation which we claim is invertible.
  By \thmref{thm:nil-localization} it will suffice to show this after pre-composition with $v^* : \MOTSPC^\cl(S_\cl) \to \MOTSPC(S)$.
  This is immediate from the fact that $v^*$ is fully faithful (\thmref{thm:nil-localization}) and the commutativity of the right-hand square in \eqref{eq:v w squares}.
\end{proof}

\begin{rem}\label{rem:functoriality of nil-localization}
  The discussion of \ssecref{ssec:nil/functoriality} makes sense in the classical setting and provides $\MotSpc^\cl(S_\cl)$ with the same functorialities as $S$ varies.
  The equivalence of \thmref{thm:intro/comp1} is compatible with all the operations $f_\sharp$, $f^*$, $f_*$, as well as $\times$, and $\uHom$:

  \begin{defnlist}
    \item\label{item:apdfjs}
    Let $f : T \to S$ be a morphism of affine spectral schemes.
    Then we have commutative squares
    \[
      \begin{tikzcd}[column sep=5em]
        \MotSpc(S) \ar{r}{\L w_!}\ar{d}{f^*}
          & \MotSpc^\cl(S_\cl) \ar{d}{f_\cl^*}
        \\
        \MotSpc(T) \ar{r}{\L w_!}
          & \MotSpc^\cl(T_\cl),
      \end{tikzcd}
      \qquad
      \begin{tikzcd}[column sep=5em]
        \MotSpc^\cl(T_\cl) \ar{r}{w^*}\ar{d}{f_{\cl,_*}}
          & \MotSpc(T) \ar{d}{f_*}
        \\
        \MotSpc^\cl(S_\cl) \ar{r}{w^*}
          & \MotSpc(S).
      \end{tikzcd}
    \]
    Indeed, the left-hand square is induced by the commutative square
      \begin{equation*}
        \begin{tikzcd}
          \sA_{/S} \ar{r}{w}\ar{d}
            & \sA^\cl_{/S_\cl} \ar{d}
          \\
          \sA_{/T} \ar{r}{w}
            & \sA^\cl_{/T_\cl},
        \end{tikzcd}
      \end{equation*}
    where the upper horizontal arrow is (derived) base change along $f$, and the lower horizontal arrow is classical base change along $f_\cl$.
    The right-hand square is obtained by passage to right adjoints.
    
    \item
    Since the horizontal arrows in the squares above are equivalences (\thmref{thm:intro/comp1}), the squares are horizontally right- and left-adjointable, respectively.
    In other words, they give rise to further commutative squares
    \[
      \begin{tikzcd}[column sep=5em]
        \MotSpc^\cl(S_\cl) \ar{r}{w^*}\ar{d}{f_{\cl}^*}
          & \MotSpc(S) \ar{d}{f^*}
        \\
        \MotSpc^\cl(T_\cl) \ar{r}{w^*}
          & \MotSpc(T),
      \end{tikzcd}
      \qquad
      \begin{tikzcd}[column sep=5em]
        \MotSpc(T) \ar{r}{\L w_!}\ar{d}{f_*}
          & \MotSpc^\cl(T_\cl) \ar{d}{f_{\cl,*}}
        \\
        \MotSpc(S) \ar{r}{\L w_!}
          & \MotSpc^\cl(S_\cl).
      \end{tikzcd}
    \]

    \item    
    Similarly $f$ is smooth, then we have commutative squares
    \[
      \begin{tikzcd}[column sep=5em]
        \MotSpc(T) \ar{r}{\L w_!}\ar{d}{f_\sharp}
          & \MotSpc^\cl(T_\cl) \ar{d}{f_{\cl,\sharp}}
        \\
        \MotSpc(S) \ar{r}{\L w_!}
          & \MotSpc^\cl(S_\cl),
      \end{tikzcd}
      \qquad
      \begin{tikzcd}[column sep=5em]
        \MotSpc^\cl(T_\cl) \ar{r}{w^*}\ar{d}{f_{\cl,\sharp}}
          & \MotSpc(T) \ar{d}{f_\sharp}
        \\
        \MotSpc^\cl(S_\cl) \ar{r}{w^*}
          & \MotSpc(S).
      \end{tikzcd}
    \]
    Here the left-hand square is induced by the commutative square
    \begin{equation*}
      \begin{tikzcd}
        \sA_{/T} \ar{r}{v}\ar{d}
          & \sA^\cl_{/T_\cl} \ar{d}
        \\
        \sA_{/S}\ar{r}{v}
          & \sA^\cl_{/S_\cl}.
      \end{tikzcd}
    \end{equation*}
    The right-hand square comes from the horizontal right-adjointability of the left-hand one.

    \item\label{item:aufpasdmf}
    Finally, consider the operations $\times$ and $\uHom$.
    Note that $w^*$ is cartesian monoidal since it is a right adjoint.
    Its left adjoint $\L w_!$ is also monoidal, since $w_!$ is clearly monoidal (as can be checked on representables) and $\L$ preserves finite products by \remref{rem:SPC localization basics}(iii).
    Then by adjunction we have a canonical isomorphism
    \[
      w^* \big(\uHom(\L w_!(\sF), \sG)\big) \to \uHom(\sF, v^*(\sG))
    \]
    in $\MotSpc(S)$, for any $\sF \in \MotSpc(S)$ and $\sG \in \MotSpc^\cl(S_\cl)$.
    This induces in turn for every $\sF,\sG\in\MotSpc(S)$ isomorphisms
    \begin{align*}
      \uHom(\L w_!(\sF), \L w_!(\sG))
      &\simeq \L w_!w^* \big(\uHom(\L w_!(\sF), \L w_!(\sG))\big)\\
      &\simeq \L w_! \big(\uHom(\sF, v^*\L w_!(\sG))\big)\\
      &\simeq \L w_! \big(\uHom(\sF, \sG)\big)
    \end{align*}
    where the first and third isomorphisms come from \thmref{thm:intro/comp1}.
  \end{defnlist}
\end{rem}


\ssec{\texorpdfstring{$\bV$}{V}-linear motivic objects}
\label{ssec:nil/V-linear}

Let $S$ be an affine spectral scheme and let $\sC_{/S} \sub \Aff_{/S}$ be an admissible subcategory.
By replacing the \inftyCat of spaces in \defnref{defn:A-fibred motivic space} with any given presentable \inftyCat $\bV$, we can define a $\bV$-linear variant of the construction $\MotSpc(\sC_{/S})$:

\begin{defn}\label{defn:V-linear categories}
  A \emph{$\sC$-fibred motivic $\bV$-object} is a $\bV$-valued presheaf $(\sC_{/S})^\op \to \bV$ satisfying $\A^1$-homotopy invariance and Nisnevich excision.
  We write $\MotSpc(\sC_{/S})_\bV$ for the \inftyCat of motivic $\sC$-fibred $\bV$-objects.
\end{defn}

Similarly, given an admissible subcategory $\sC^\cl_{/S_\cl} \sub \AffCl_{/S_\cl}$, we may consider the \inftyCat $\MotSpc^\cl(S_\cl)_\bV$ of $\sC^\cl$-fibred motivic $\bV$-objects.
We have $\bV$-linear analogues of each of the categories defined in \emph{loc. cit.}:
\[
  \Spc(\sC_{/S})_\bV, \quad \Spc(\sC^\cl_{/S_\cl})_\bV,
  \quad \MotSpc(\sC_{/S})_\bV, \quad \MotSpc^\cl(\sC^\cl_{/S_\cl})_\bV.
\]

\begin{exam}\label{exam:motivic S^1-spectra}
  Taking $\bV$ to be the stable presentable \inftyCat $\Spt$ of spectra, we obtain \inftyCats of fibred motivic spectra.
  We will refer to these as fibred \emph{motivic $S^1$-spectra}, to distinguish them from the notion of motivic spectra with respect to the Thom space of the trivial line bundle.
\end{exam}

\begin{rem}\label{rem:V stable}
  If $\bV$ is \emph{stable}, then so is $\MotSpc(\sC_{/S})_\bV$ for any admissible subcategory $\sC_{/S} \sub \Aff_{/S}$.
\end{rem}

\begin{rem}\label{rem:tensor V}
  The \inftyCat $\MotSpc(\sC_{/S})_\bV$ can also be described as the tensor product of presentable \inftyCats $\MotSpc(\sC_{/S}) \otimes \bV$ in the sense of \cite[Sect.~4.8]{HA-20170918}.
  An analogous description holds for the classical variant $\MotSpc(\sC^\cl_{/S_\cl})$, for $\sC^\cl_{/S_\cl} \sub \AffCl_{/S_\cl}$ admissible.
  It follows that when $\bV$ is presentably symmetric monoidal, these categories also inherit presentably symmetric monoidal structures.
  \remref{rem:i^*Hom} then carries over to the $\bV$-linear setting.
\end{rem}

The description of \remref{rem:tensor V} immediately gives the following generalization of the comparison result of \thmref{thm:intro/comp1} (or rather the more precise statement proven in \ssecref{ssec:nil/comparison}):

\begin{thm}\label{thm:nil-localization V-linear}
  Let $S$ be a spectral affine scheme and let $\sA_{/S} \sub \Aff_{/S}$ be a \emph{narrow} subcategory.
  Let $\sA^\cl_{/S_\cl} \sub \AffCl_{/S}$ be as in \examref{exam:A^cl} and let $w : \sA_{/S} \to \sA_{/S_\cl}$ denote the classical truncation functor (\constrref{constr:v}).
  Then for any presentable \inftyCat $\bV$, the adjunction
  \begin{equation*}
    \L w_! : \MotSpc(S)_\bV \to \MotSpc^\cl(S_\cl)_\bV,
    \quad
    w^* : \MotSpc^\cl(S_\cl)_\bV \to \MotSpc(S)_\bV
  \end{equation*}
  is an equivalence of \inftyCats.
  Moreover, this equivalence is compatible with the operations $f^*$, $f_*$ for any morphism $f : T \to S$, with $f_\sharp$ when $f$ exhibits $T$ as an object of $\sA_{/S}$, with products and with internal homs.
\end{thm}

Finally, let us note the following two properties which are specific to the stable case.

\begin{prop}\label{prop:Lmot of spectra}
  Let $\bV$ be a \emph{stable} presentable \inftyCat.
  Let $\sC_{/S} \sub \Aff_{/S}$ be an admissible subcategory.
  Then on the \inftyCat of $\sC$-fibred $\bV$-objects, the $\A^1$-localization functor $\Lhtp$ preserves the property of Nisnevich excision.
  In particular, the motivic localization functor $\L$ can be computed by the formula
  \begin{equation*}
    \L \simeq \Lhtp \LNis.
  \end{equation*}
\end{prop}

\begin{proof}
  Let $\sF$ be a $\sC$-fibred $\bV$-object.
  If $\sF$ is Nisnevich-excisive, then its $\A^1$-localization $\Lhtp(\sF)$ is still Nisnevich-excisive, in view of the formula \eqref{eq:formula for Lhtp} and the fact that colimits commute with finite limits in stable \inftyCats.
  Therefore the claim follows from \remref{rem:SPC localization basics} (iii).
\end{proof}

\begin{cor}\label{cor:Lhtp of spectra using v_!}
  Let $\bV$ be a \emph{stable} presentable \inftyCat.
  Let $\sF \in \MotSpc(S)_\bV$ be an $\sA$-fibred $\bV$-object over $S$ (where $\sA_{/S} \sub \Aff_{/S}$ is narrow).
  If $\sF$ is Nisnevich-excisive, then the canonical map
  \begin{equation*}
    \Lhtp(\sF) \to w^*\Lhtp w_!(\sF)
  \end{equation*}
  is invertible.
\end{cor}

\begin{proof}
  It follows from \thmref{thm:nil-localization V-linear} (cf. \corref{cor:Lmot(F) in terms of Lmot v_!(F)}) that the canonical map $\L(\sF) \to w^*\L w_!(\sF)$ is invertible.
  Since $\sF$ and hence $w_!(\sF)$ are Nisnevich-excisive, we conclude by \propref{prop:Lmot of spectra}.
\end{proof}


\section{Comparison with \texorpdfstring{$\A^{1,\flat}$}{A1flat}-motivic homotopy theory}
\label{sec:smflat}

\ssec{Flat affine spaces}

\begin{notat}
  Let $R$ be an \Ering.
  Denote by $R[T_1,\ldots,T_n]$ denote the polynomial $R$-algebra in $n$ variables $T_1,\ldots,T_n$ (in degree zero).
  This is by definition the monoid algebra $R[\bN^n] = R \otimes \Sigma^\infty_+(\bN^n)$, where $\bN$ is the set of natural numbers, viewed as a discrete (additive) $\Einfty$-monoid space.
  Note that we have canonical isomorphisms $\pi_*(R[T_1,\ldots,T_n]) \simeq \pi_*(R) \otimes_{\pi_0(R)} \pi_0(R)[T_1,\ldots,T_n]$, so that $R[T_1,\ldots,T_n]$ is flat over $R$.
\end{notat}

\begin{defn}\label{defn:flat affine space}
  For every $n\ge 0$, let $\A^{n,\flat}$ denote the affine spectral scheme $\Spec(\bS[T_1,\ldots,T_n])$, where $\bS$ is the sphere spectrum.
  Note that we have a canonical isomorphism $(\A^{n,\flat})_\cl \simeq \A^n_\cl$.
  We refer to $\A^{n,\flat}$ as the \emph{flat affine space} (over the sphere spectrum).
  If $S$ is classical, then $S\times\A^{n,\flat} = S\times\A^n_{\cl}$.
\end{defn}

\begin{rem}\label{rem:A^n_flat structure}
  The affine spectral schemes $\A^{n,\flat}$ are equipped with the following additional structure:
  \begin{defnlist}
    \item
    The flat affine line $\A^{1,\flat}$ has the structure of a commutative monoid under the operation of multiplication.
    This is induced by the cocommutative comonoid structure on the commutative monoid $\bN$.
    For example, the multiplication morphism $\A^{1,\flat} \times \A^{1,\flat} \to \A^{1,\flat}$ corresponds to the diagonal $\Sigma^\infty_+(\bN) \to \Sigma^\infty_+(\bN\times\bN) \simeq \Sigma^\infty_+(\bN) \otimes \Sigma^\infty_+(\bN)$.
    Similarly the counit of $\bN$ induces an $\Einfty$-ring homomorphism $\bS[T] \simeq \Sigma^\infty_+(\bN) \to \Sigma^\infty_+(\pt) \simeq \bS$ which corresponds to the unit section $s_1 : \Spec(\bS) \to \A^{1,\flat}$.

    \item
    The flat affine line $\A^{1,\flat}$ also admits a zero section $s_0 : \Spec(\bS) \to \A^{1,\flat}$, which can be constructed as follows.
    Identify the discrete pointed $\Einfty$-monoid space $\pt_+$ with the set $\{0, 1\}$, viewed as a multiplicative monoid with base point $0$ and identity element $1$.
    Since $\bN$ is freely generated as a (discrete) commutative monoid by the element $1 \in \bN$, either choice of element $i \in \{0,1\}$ gives rise to a unique homomorphism $\sigma'_i : \bN_+ \to \pt_+$ of pointed commutative monoids such that $\sigma'_i(1) = i$.
    Regarding $\sigma'_i$ as a homomorphism of discrete pointed $\Einfty$-monoid spaces, application of the symmetric monoidal functor $\Sigma^\infty$ produces $\Einfty$-ring homomorphisms
      \begin{equation*}
        \sigma_i : \bS[T] \simeq \Sigma^\infty(\bN_+) \to \Sigma^\infty(\pt_+) \simeq \bS
      \end{equation*}
    for each $i \in \{0,1\}$.
    For $i=1$ this is the same homomorphism defining the unit section $s_1$, and we let $s_0$ denote the section $\Spec(\bS) \to \Spec(\bS[T]) = \A^{1,\flat}$ corresponding to $\sigma_0$.

    \item
    For every $n\ge 0$, the flat affine space $\A^{n,\flat}$ admits the structure of a module over $\A^{1,\flat}$.
    This is induced by the canonical comodule structure on the commutative monoid $\bN^n$ over the comonoid $\bN$, where the coaction homomorphism sends $\bN^n \to \bN \times \bN^n$
    \[
      (k_1,k_2,\ldots,k_n) \mapsto (k_1+\cdots+k_n, k_1,k_2,\ldots,k_n).
    \]
    In particular, there is a ``scalar multiplication'' morphism
    \[
      \A^{1,\flat} \times \A^{n,\flat} \to \A^{n,\flat}.
    \]

    \item\label{item:lzcxhvkjz}
    After base change along $\Spec(\bZ) \to \Spec(\bS)$, the flat affine spaces $\Spec(\bZ) \times \A^{n,\flat} \simeq \A^n_{\cl}$ become abelian groups under addition.
  \end{defnlist}
\end{rem}

\begin{rem}\label{rem:A^1_flat interval}
  The zero section of $\A^{1,\flat}$ is compatible with the multiplicative structure, in the sense that $\A^{1,\flat}$ defines an \emph{interval object} in the sense of Morel--Voevodsky \cite[Sect.~2.3]{MorelVoevodsky}.
\end{rem}

\ssec{\texorpdfstring{$\Sm^\flat$}{Smflat}-fibred motivic spaces}

Unlike $S \times \A^n$, the flat affine spaces $S\times\A^{n,\flat}$ are not smooth over $S$ in the sense of \defnref{defn:smooth} (except under the conditions of \remref{rem:A^n_flat to A^n}).
However, they are smooth in the following sense:

\begin{defn}\label{defn:fibre-smooth}
A morphism of affine spectral schemes $X \to S$ is called \emph{fibre-smooth} if it is almost of finite presentation, flat, and on classical truncations induces a morphism $X_\cl \to S_\cl$ that is smooth in the sense of classical algebraic geometry.
\end{defn}

\begin{rem}
  From \cite[Cor.~11.2.4.2]{SAG-20180204} and \cite[\sectsign~17.3]{EGAIV4} it follows that a morphism of affine spectral schemes $\Spec(B) \to \Spec(A)$ is fibre-smooth if and only if $A \to B$ is fibre-smooth in the sense of \cite[Def.~11.2.3.1]{SAG-20180204}.
\end{rem}

\begin{rem}
  Let $f : X \to S$ be a fibre-smooth morphism.
  Then Zariski-locally on $X$, there exists an étale $S$-morphism $X \to S \times \A^{n,\flat}$ for some $n\ge 0$.
  This follows from Remark~11.2.3.5 and Proposition~11.2.4.1 of \cite{SAG-20180204}, combined with \cite[Thm.~7.5.0.6]{HA-20170918}.
  Contrast with \remref{rem:smooth factor}.
\end{rem}

The following is the same as \defnref{defn:adm} except for the last condition.

\begin{defn}\label{defn:flatadm}
  We say that a full subcategory $\sC^\flat_{/S} \sub \Aff_{/S}$ is \emph{$\flat$-admissible} if it is essentially small and satisfies the following conditions:
  \begin{defnlist}[series=adm]
    \item\label{item:flatadm1}
    The affine spectral scheme $S$ (viewed over $S$ via the identity) belongs to $\sC^\flat_{/S}$.
    \item\label{item:flatadm2}
    If $X$ belongs to $\sC^\flat_{/S}$ and $Y$ is étale over $X$, then $Y$ belongs to $\sC^\flat_{/S}$.
    \item\label{item:flatadm3}
    If $X$ belongs to $\sC^\flat_{/S}$, then $X \times \A^{n,\flat}$ belongs to $\sC^\flat_{/S}$ for every $n\ge 0$.
  \end{defnlist}

  A \emph{$\flat$-narrow} subcategory $\sC^\flat_{/S} \sub \Aff_{/S}$ is a $\flat$-admissible subcategory which is contained in the full subcategory $\Sm^\flat_{/S} \sub \Aff_{/S}$ of fibre-smooth spectral affine schemes over $S$.
\end{defn}

\begin{exam}\label{exam:minimal flatadm}
  Let $\sA^{\flat, 0}_{/S} \sub \Aff_{/S}$ denote the full subcategory spanned by $X \in \Aff_{/S}$ for which the structural morphism $X \to S$ factors through an étale morphism
  \[ X \to S \times \A^{n,\flat} \]
  over $S$.
  Then $\sA^{\flat,0}_{/S}$ is the minimal $\flat$-admissible subcategory of $\Aff_{/S}$ (same proof as \examref{exam:minimal adm}).
\end{exam}

\begin{rem}\label{rem:C^cl from C^flat}
  Let $S$ be a spectral affine scheme.
  Let $\sC^\flat \sub \Aff_{/S}$ be a $\flat$-admissible subcategory.
  As in \examref{exam:A^cl}, the full subcategory $\sC^\cl_{/S_\cl} \sub \AffCl_{/S_\cl}$ spanned by classical truncations of objects in $\sC^\flat$ is admissible.
\end{rem}

\begin{defn}\label{defn:C^flat fibred}
  Let $S$ be a spectral affine scheme and $\sC^\flat_{/S} \sub \Aff_{/S}$ a $\flat$-admissible subcategory.
  Let $\sF$ be a $\sC^\flat$-fibred space (\defnref{defn:fibred space}), i.e., a presheaf of spaces on $\sC^\flat_{/S}$.
  We say that $\sF$ is Nisnevich excisive if it is reduced and sends Nisnevich squares $Q$ to cartesian squares $\Gamma(Q, \sF)$ (\defnref{defn:Nisnevich}).
  We say that $\sF$ satisfies \emph{$\A^{1,\flat}$-homotopy invariance} if for any $X \in \sC^\flat_{/S}$, the canonical map $\Gamma(X, \sF) \to \Gamma(X\times\A^{1,\flat},\sF)$ is invertible, where $\A^{1,\flat}$ denotes the flat affine line (\defnref{defn:flat affine space}).
  We say that $\sF$ is \emph{motivic} if it is Nisnevich excisive and $\A^{1,\flat}$-homotopy invariant.
  We denote by $\MotSpc(\sC^\flat_{/S}) \subseteq \Spc(\sC^\flat_{/S})$ the full subcategory of motivic $\sC^\flat$-fibred spaces.
\end{defn}

\begin{notat}\label{notat:ambig}
  If $\sC_{/S} \sub \Aff_{/S}$ is both admissible and $\flat$-admissible, then there is possible ambiguity in the terminology ``$\sC$-fibred motivic space'' and in the notation $\MotSpc(\sC_{/S})$.
  To maintain the distinction we introduce the following convention: we write $\sC^\flat_{/S} \sub \Aff_{/S}$ for the same full subcategory when it is to be regarded as a $\flat$-admissible subcategory.
  Thus a \emph{motivic $\sC$-fibred space} is an $\A^1$-invariant Nisnevich-excisive $\sC$-fibred space, while a \emph{motivic $\sC^\flat$-fibred space} is an $\A^{1,\flat}$-invariant Nisnevich-excisive $\sC$-fibred space.
  In particular, $\MotSpc(\sC_{/S})$ and $\MotSpc(\sC^\flat_{/S})$ are two distinct full subcategories of $\Spc(\sC_{/S}) = \Spc(\sC^{\flat}_{/S})$.
\end{notat}

\begin{rem}\label{rem:SPC flat localization}
  For any $\flat$-admissible subcategory $\sC^\flat_{/S} \sub \Aff_{/S}$, the full subcategories of Nisnevich-excisive, $\A^{1,\flat}$-homotopy invariant, and motivic $\sC^\flat$-fibred spaces are each left Bousfield localizations of the \inftyCat of $\sC^\flat$-fibred spaces.
  We write $\LNis$, $\Lhtpflat$, and $\L^\flat$ for the respective localization functors.
  The remarks in \ref{rem:SPC localization basics} remain valid \emph{mutatis mutandis} up to replacing the spectral affine spaces $\A^n$ by the flat affine spaces $\A^{n,\flat}$ in item~(ii).
\end{rem}

\begin{prop}\label{prop:MotSpc^flat vs MOTSPC}
  Let $S$ be a spectral affine scheme.
  Fix an inclusion $\sC^\flat_{/S} \sub \sD^\flat_{/S}$ of $\flat$-admissible subcategories of $\Aff_{/S}$.
  Denote by $\iota : \sC^\flat_{/S} \hook \sD^\flat_{/S}$ the inclusion functor and by $\iota_! : \Spc(\sC^\flat_{/S}) \to \Spc(\sD^\flat_{/S})$ its left Kan extension.
  Then the assignment $\sF \mapsto \L^\flat\iota_!(\sF)$ induces a fully faithful functor of \inftyCats
  \begin{equation*}
    \L^\flat\iota_! : \MotSpc(\sC^\flat_{/S}) \to \MotSpc(\sD^\flat_{/S}),
  \end{equation*}
  whose essential image is generated under sifted colimits by objects of the form $\L^\flat (\h_S(X))$, where $X$ belongs to $\sC_{/S}$.
\end{prop}
\begin{proof}
  Same proof as \propref{prop:MotSpc vs MOTSPC}.
\end{proof}

\begin{rem}
  Let $\sC^\flat_{/S} \sub \Aff_{/S}$ be $\flat$-admissible.
  Using the interval structure on $\A^{1,\flat}$ (Remarks~\ref{rem:A^n_flat structure} and \ref{rem:A^1_flat interval}), we can make sense of \emph{$\A^{1,\flat}$-homotopies} between morphisms of $\sC^\flat$-fibred spaces, and therefore of \emph{strict $\A^{1,\flat}$-homotopy equivalences} (see \cite[Def.~2.3.6]{KhanLocalization}).
\end{rem}

\ssec{\texorpdfstring{$\A^1$}{A1}-homotopies vs. \texorpdfstring{$\A^{1,\flat}$}{A1flat}-homotopies}

\begin{rem}\label{rem:A^n_flat to A^n}
For any \cnEring $R$ and integer $n\ge 0$, there is a canonical homomorphism of $\Einfty$-$R$-algebras $R\{T_1,\ldots,T_n\} \to R[T_1,\ldots,T_n]$ determined by the assignment $T_i \mapsto T_i$.
This gives rise to canonical morphisms
  \begin{equation*}
    \varepsilon_S : S\times \A^{n,\flat} \to S\times \A^n
  \end{equation*}
for every affine spectral scheme $S$, which are invertible if and only if either $n=0$ or $S$ is of characteristic zero (i.e., $S = \Spec(R)$ with $R$ an $\Einfty$-$\bQ$-algebra).
\end{rem}

\begin{rem}\label{rem:A^1_flat to A^1 intervals}
  The map $\varepsilon : \A^{1,\flat} \to \A^1$ is compatible with interval structures.
  That is, it preserves the zero and unit sections, and is compatible with the multiplicative structures in the sense that the diagram
  \begin{equation*}
    \begin{tikzcd}
      \A^{1,\flat} \times \A^{1,\flat} \ar{r}\ar{d}
        & \A^{1,\flat} \ar{d}
      \\
      \A^1 \times \A^1 \ar{r}
        & \A^1
    \end{tikzcd}
  \end{equation*}
  commutes.
  By restriction of scalars, we may therefore regard $\A^n$ as a module over $\A^{1,\flat}$.
  Note that for every $n\geq 0$, the morphism $\varepsilon : \A^{n,\flat} \to \A^n$ is then $\A^{1,\flat}$-linear.

  Note also that over $\Spec(\bZ)$, $\varepsilon$ defines a group homomorphism
  \[
    \A^n_{\cl} \simeq \Spec(\bZ) \times \A^{n,\flat} \to \Spec(\bZ) \times \A^n
  \]
  with respect to the additive structures (\remref{rem:A^n_flat structure}\itemref{item:lzcxhvkjz}).
\end{rem}

\begin{lem}\label{lem:A^1 vs A^1_flat-invariance}
  Let $S$ be an affine spectral scheme.
  Let $\sC_{/S}$ be an admissible and $\flat$-admissible subcategory of $\Aff_{/S}$.
  Then we have:
  \begin{thmlist}
    \item
    Every $\A^1$-local equivalence between $\sC$-fibred spaces over $S$ is an $\A^{1,\flat}$-local equivalence.

    \item
    Every $\A^{1,\flat}$-homotopy invariant $\sC$-fibred space over $S$ is $\A^1$-homotopy invariant.
  \end{thmlist}
  In particular, there is an inclusion
  \[ \MotSpc(\sC^\flat_{/S}) \sub \MotSpc(\sC_{/S}) \]
  of subcategories of $\Spc(\sC_{/S})$.
\end{lem}

\begin{proof}
  By adjunction, the two statements are equivalent.
  To prove the first it will suffice to show that for every $X \in \sC_{/S}$, the morphism $\h_S(X \times \A^1) \to \h_S(X)$ is an $\A^{1,\flat}$-local equivalence.
  In fact, we claim that every strict $\A^1$-homotopy equivalence is a strict $\A^{1,\flat}$-homotopy equivalence (hence \emph{a fortiori} an $\A^{1,\flat}$-local equivalence).
  Indeed, the canonical map $\varepsilon_S : S\times\A^{1,\flat} \to S\times\A^1$ (\remref{rem:A^n_flat to A^n}) is a morphism of interval objects (\remref{rem:A^1_flat to A^1 intervals}), so composition with $\varepsilon_S$ sends elementary $\A^1$-homotopies to elementary $\A^{1,\flat}$-homotopies.
\end{proof}

The following key lemma shows that the comparison morphism $\varepsilon : \A^{n,\flat} \to \A^n$ is a ``universal'' $\A^{1,\flat}$-local equivalence, at least over $\Spec(\bZ)$.
The reason for this restriction is that the proof uses the additive structure on the flat affine spaces (\remref{rem:A^n_flat structure}\itemref{item:lzcxhvkjz}).

\begin{lem}\label{lem:epsilon A^1_flat-equiv}
  Let $S$ be a spectral affine scheme defined over $\Spec(\bZ)$.
  Let $\sC_{/S} \sub \Aff_{/S}$ be an admissible and $\flat$-admissible subcategory.
  Given $X\in\sC_{/S}$ and an $S$-morphism $f : X \to \A^n_S$, consider the cartesian square
  \[
    \begin{tikzcd}
      X^\flat \ar{r}\ar{d}{g}
        & X \ar{d}{f}
      \\
      S \times \A^{n,\flat} \ar{r}{\varepsilon_S}
        & S \times \A^n.
    \end{tikzcd}
  \]
  Then the morphism $\h_S(X^\flat) \to \h_S(X)$ is an $\A^{1,\flat}$-local equivalence of $\sC$-fibred spaces.
\end{lem}
\begin{proof}
  Note that $f$ induces a section $s : X \to (S\times \A^n)\fibprod_S X = X \times \A^n$ of the projection $X \times \A^n \to X$, and similarly $g$ induces a section $t : X \to X\times \A^{n,\flat}$.
  These fit into a factorization of the given square:
  \[
    \begin{tikzcd}
      X^\flat \ar{r}\ar{d}{t}
        & X \ar{d}{s}
      \\
      X \times \A^{n,\flat} \ar{r}{\varepsilon_X}\ar{d}
        & X\times \A^n\ar{d}
      \\
      S \times \A^{n,\flat} \ar{r}{\varepsilon_S}
        & S\times \A^n
    \end{tikzcd}
  \]
  Since the lower square is cartesian, so is the upper square.
  Therefore we may as well assume that $X = S$, in which case the claim becomes that $\h_X(X^\flat)$ is $\A^{1,\flat}$-contractible.
  Since $S$ is defined over $\Spec(\bZ)$, $\varepsilon$ is a group homomorphism with respect to the additive structures (\remref{rem:A^1_flat to A^1 intervals}).
  Therefore it gives rise to an isomorphism between $X^\flat$ and the fibre of $\varepsilon$ over the zero section.
  Thus we may also assume that $s$ is the zero section.

  Since $s : X \to \A^n_X$ lifts to the zero section $s : X \to \A^n_{\flat,X}$, it induces an $X$-morphism $X \to X^\flat$ and hence a base point of the $\sC$-fibred space $\h_X(X^\flat)$.
  It will suffice to exhibit an $\A^{1,\flat}$-homotopy contracting $\h_X(X^\flat)$ to this base point.
  Recall that $\A^{1,\flat}$ acts compatibly on $\A^n_{\flat,X}$ and $\A^n_X$ (Remarks~\ref{rem:A^n_flat structure} and \ref{rem:A^1_flat to A^1 intervals}).
  The action on the latter restricts along the zero section $s : X \to \A^n_X$ to the trivial action on $X$.
  The induced $\A^{1,\flat}$-action on $X^\flat$ is a morphism
  \[
    \A^{1,\flat} \times X^\flat \to X^\flat,
  \]
  which induces the $\A^{1,\flat}$-homotopy desired.
\end{proof}

\begin{cor}\label{cor:image of Hflat}
  Let $S$ be a spectral affine scheme over $\Spec(\bZ)$.
  Let $\sA^\flat_{/S} \sub \Aff_{/S}$ be a $\flat$-narrow subcategory and let $\sC_{/S} \sub \Aff_{/S}$ be an admissible and $\flat$-admissible subcategory containing $\sA^\flat_{/S}$.
  Let $\iota : \sA^\flat_{/S} \hook \sC_{/S}$ denote the inclusion.
  Then the essential image of the fully faithful functor (\propref{prop:MotSpc^flat vs MOTSPC})
  \begin{equation*}
    \Lmotflat\iota_! : \MotSpc(\sA^\flat_{/S}) \to \MotSpc(\sC^\flat_{/S})
  \end{equation*}
  is generated under sifted colimits by objects of the form $\Lmotflat \h_S(X)$, where $X \in \Sm_{/S}$ is affine and admits an étale $S$-morphism to the \emph{spectral} affine space $\A^n_{S}$, for some $n\ge 0$.
\end{cor}

\begin{proof}
  Let $\sE$ and $\sE_\flat$ denote the full subcategories of $\MotSpc(\sC^\flat_{/S})$ generated under sifted colimits by objects of the form $\Lmotflat \h_S(X)$, where $X$ admits an étale $S$-morphism to $S\times\A^n$, respectively to $S\times\A^{n,\flat}$, for some $n\ge 0$.
  We may as well assume $\sA^\flat_{/S}$ is the minimal $\flat$-admissible subcategory (\examref{exam:minimal flatadm}), spanned by $X \in \Aff_{/S}$ which admit an étale morphism $X \to S \times \A^{n,\flat}$ for some $n\ge 0$.
  Then by \propref{prop:MotSpc^flat vs MOTSPC}, $\sE_\flat$ is identified with the essential image of the functor in question, so it will suffice to show that $\sE = \sE_\flat$.

  Let $X \in \Aff_{/S}$.
  Suppose $X$ admits an étale $S$-morphism $f : X \to S\times\A^n$ for some $n\ge 0$.
  The base change of $f$ along $\varepsilon_S : S\times \A^{n,\flat} \to S\times\A^n$ is an étale morphism $X^\flat \to S\times \A^{n,\flat}$.
  By \lemref{lem:epsilon A^1_flat-equiv}, the canonical morphism $X^\flat \to X$ induces an isomorphism $\Lmotflat \h_S(X^\flat)\simeq\Lmotflat \h_S(X)$, hence in particular $\Lmotflat \h_S(X) \in \sE_\flat$.
  This shows $\sE \subseteq \sE_\flat$.

  For the other direction, suppose $X \in \Aff_{/S}$ admits an étale $S$-morphism $g : X \to S\times \A^{n,\flat}$ for some $n\ge 0$.
  Since $\varepsilon_S$ is an isomorphism on classical truncations, it follows from \cite[Thm.~7.5.0.6]{HA-20170918} that there exists $Y\in\Aff_{/S}$, an étale $S$-morphism $f : Y \to S\times\A^n$, and a cartesian square
  \[
    \begin{tikzcd}
      X \ar{r}\ar{d}{g}
        & Y \ar{d}{f}
      \\
      S\times \A^{n,\flat} \ar{r}{\varepsilon_S}
        & S\times\A^n.
    \end{tikzcd}
  \]
  By \lemref{lem:epsilon A^1_flat-equiv}, there is an isomorphism $\Lmotflat \h_S(X) \simeq \Lmotflat \h_S(Y)$, hence in particular $\Lmotflat \h_S(X) \in \sE$.
  This shows $\sE_\flat \subseteq \sE$.
\end{proof}

\ssec{The comparison}
\label{ssec:smflat/comp}

In this subsection we prove the following statement, which in particular yields the equivalence between (ii) and (iii) in \thmref{thm:intro/comp1}.

\begin{thm}\label{thm:comp2}
  Let $S$ be a spectral affine scheme defined over $\Spec(\bZ)$.
  Let $\sA^\flat_{/S}$ be any $\flat$-narrow subcategory of $\Aff_{/S}$ and let $w^\flat : \sA^\flat_{/S} \to \sA^\cl_{/S}$ be the restriction of the classical truncation functor $v : \sB_{/S} \to \sB^\cl_{/S_\cl}$.
  Then the adjunction
  \[
    \L w^\flat_! : \MotSpc(\sA^\flat_{/S}) \to \MotSpc(\sA^\cl_{/S})
    \quad
    w^* : \MotSpc(\sA^\cl_{/S}) \to \MotSpc(\sA^\flat_{/S})
  \]
  is an equivalence.
\end{thm}

We fix the following notation for this subsection.

\begin{notat}\label{notat:pasidfop}
  Let $S$ be a spectral affine scheme.
  Let $\sA^\flat_{/S} \sub \Aff_{/S}$ be a $\flat$-narrow subcategory and $\sB \sub \Aff_{/S}$ be a $\flat$-admissible and broad subcategory containing $\sA^\flat_{/S}$.
  Let $\iota : \sA^\flat_{/S} \hook \sB_{/S}$ denote the inclusion.
  Let $\sB^\flat_{/S}$ be as in \remref{notat:ambig} and let $\sB^\cl_{/S} \sub \AffCl_{/S}$ be as in \examref{exam:A^cl}.
  We write
  \[
    \MotSpc(\sB_{/S}),
    \quad \MotSpc(\sB^\flat_{/S}),
    \quad \MotSpc(\sB^\cl_{/S_\cl})
  \]
  for the \inftyCats of motivic spaces formed respectively out of the admissible subcategory $\sB_{/S} \sub \Aff_{/S}$, the $\flat$-admissible subcategory $\sB^\flat_{/S} \sub \Aff_{/S}$, and the admissible subcategory $\sB^\cl_{/S_\cl} \sub \AffCl_{/S}$.
\end{notat}

\begin{rem}\label{rem:v_!,v^* flat}
  By construction the classical truncation functor $v : \Aff_{/S} \to \AffCl_{/S_\cl}$ (\constrref{constr:v}) restricts to a functor $v : \sB_{/S} \to \sB^\cl_{/S_\cl}$.
  Since the latter preserves Nisnevich squares and sends $\A^{1,\flat}$ to $\A^1_\cl$, one sees as in \lemref{lem:v_!,v^* and motivic equivalences} that $v$ induces an adjunction
  \[
    \L v_! : \MotSpc(\sB^\flat_{/S}) \to \MotSpc(\sB^\cl_{/S_\cl}),
    \quad
    v^* : \MotSpc(\sB^\cl_{/S_\cl}) \to \MotSpc(\sB^\flat_{/S}).
  \]
  \propref{prop:v_!,v^* and motivic equivalences} also holds \emph{mutatis mutandis} for the functors
    \[
      v_! : \Spc(\sB^\flat_{/S}) \to \Spc(\sB^\cl_{/S_\cl}),
      \quad v^* : \Spc(\sB^\cl_{/S_\cl}) \to \Spc(\sB^\flat_{/S}).
    \]
    In particular, $v^*$ commutes with $\L^\flat$.
    By \remref{rem:v^*(F) is nil-local} this implies that $\L^\flat$ preserves nil-local objects.
\end{rem}

\begin{cor}\label{cor:A1flat Asp-fibred nil-localization}
  The functor \[ v^* : \MotSpc(\sB^\cl_{/S_\cl}) \to \MotSpc(\sB^\flat_{/S}) \] is fully faithful, with essential image spanned by the nil-local objects of $\MotSpc(\sB^\flat_{/S})$.
\end{cor}

\begin{proof}
  Follows immediately by combining \thmref{thm:nil-localization} with \lemref{lem:A^1 vs A^1_flat-invariance}.
\end{proof}

\begin{cor}\label{cor:flat nil descent}
  If $S$ is defined over $\Spec(\bZ)$, then for every $\sF \in \MotSpc(\sA^\flat_{/S})$, the $\sB^\flat$-fibred motivic space $\L^\flat \iota_!(\sF) \in \MotSpc(\sB^\flat_{/S})$ is nil-local. 
\end{cor}
\begin{proof}
  By \thmref{thm:nil descent}, every object in the essential image of \[ \L \iota_! : \MotSpc(\sA^0_{/S}) \to \MotSpc(\sB_{/S}) \] is nil-local, where $\sA^0_{/S}$ is the minimal admissible subcategory of $\Aff_{/S}$ (\examref{exam:minimal adm}).
  It follows from \remref{rem:v_!,v^* flat} that the same holds for objects in the essential image of \[ \L^\flat \L \iota_! \simeq \L^\flat \iota_! : \MotSpc(\sA^0_{/S}) \to \MotSpc(\sB^\flat_{/S}), \] where the isomorphism is due to \lemref{lem:A^1 vs A^1_flat-invariance}.
  By \corref{cor:image of Hflat} this coincides with the essential image of the functor in question.
\end{proof}

\begin{proof}[Proof of \thmref{thm:comp2}]
  We are free to choose any $\flat$-admissible and broad subcategory $\sB_{/S} \sub \Aff_{/S}$ containing $\sA^\flat_{/S}$ as in \notatref{notat:pasidfop}.
  By \propref{prop:MotSpc^flat vs MOTSPC} we have the commutative diagram
  \[ \begin{tikzcd}[column sep=5em]
    \MotSpc(\sA^\flat_{/S}) \ar{r}{\L w^\flat_{!}}\ar[swap]{d}{\Lmotflat\iota_!}
      & \MotSpc(\sA^\cl_{/S_\cl}) \ar{d}{\L\iota_!}
    \\
    \MotSpc(\sB^\flat_{/S}) \ar{r}{\L v_!}
      & \MotSpc(\sB^\cl_{/S_\cl})
  \end{tikzcd} \]
  where the vertical arrows are fully faithful.
  Write $\langle\sA^\flat_{/S}\rangle$ for the essential image of the left-hand vertical arrow, and $\langle\sA^\cl_{/S_\cl}\rangle$ for that of the right-hand vertical arrow.
  It will suffice to show that $\L v_!$ restricts to an equivalence $\langle\sA^\flat_{/S}\rangle \simeq \langle\sA^\cl_{/S_\cl}\rangle$.
  By \corref{cor:flat nil descent}, every object in the source is nil-local.
  Thus the claim follows from \corref{cor:A1flat Asp-fibred nil-localization}.
\end{proof}

\begin{cor}\label{cor:comp2 asp}
  For any quasi-compact quasi-separated spectral algebraic space $S$ over $\Spec(\bZ)$, there are canonical equivalences of \inftyCats
  \[ \MotSpc^\flat(S) \simeq \MotSpc^\cl(S_\cl). \]
  Here $\MotSpc^\flat(S)$ is the \inftyCat of $\A^{1,\flat}$-invariant Nisnevich sheaves on the site of (quasi-compact quasi-separated) fibre-smooth spectral algebraic spaces over $S$ (where fibre-smoothness is defined as in \cite[Def.~11.2.5.5]{SAG-20180204}).
\end{cor}
\begin{proof}
  Follows from \thmref{thm:comp2} by Nisnevich descent as in \corref{cor:comp asp}.
\end{proof}


\section{The Bass construction}
\label{sec:bass}


\ssec{Localizing invariants}

We briefly establish our notations and conventions for localizing invariants.
We work with $R$-linear stable \inftyCats, over a fixed \cnEring $R$, although the discussion makes sense in the greater generality of $\Perf(S)$-linear \inftyCats, for any spectral algebraic space $S$.

\begin{notat}
  Let $\Stab$ denote the \inftyCat of small stable $\infty$-categories.
  For any connective $\Einfty$-ring $R$, let $\Stab_R$ denote the \inftyCat of small stable $R$-linear $\infty$-categories.
\end{notat}

\begin{notat}
  Let $R$ be a connective $\Einfty$-ring and $E : \Stab_R \to \Spt$ a functor.
  For any $R$-algebra $A$, resp. spectral algebraic space $X$ over $\Spec(R)$, we set
  \[ E(A) := E(\Perf_A), \quad\text{resp.}~E(X) := E(\Perf(X)), \]
  where $\Perf_A$ is the stable \inftyCat of $A$-modules and $\Perf(X)$ is the stable \inftyCat of perfect complexes on $X$.
\end{notat}

\begin{rem}
Any $R$-linear stable \inftyCat $\bA$ corepresents a functor $h_\bA : \Stab_R \to \Spt$ given by the assignment
  \begin{equation*}
    h_\bA(\bA') = \Maps_{R}(\bA, \bA'),
  \end{equation*}
where $\Maps_R$ here denotes the mapping spectrum in the \inftyCat of functors $\Stab_R \to \Spt$.
\end{rem}

\begin{rem}
The Day convolution product endows the \inftyCat of functors $\Stab_R \to \Spt$ with a closed symmetric monoidal structure for which the Yoneda embedding is symmetric monoidal.
Given a functor $E : \Stab_R \to \Spt$, we write $E^\bA$ for the internal hom object $\uHom_R(h_\bA, E)$, for any $\bA \in \Stab_R$.
Note that the assignment $\bA \mapsto h_\bA$ is \emph{contravariant}, while $\bA \mapsto E^\bA$ is \emph{covariant}.
By the Yoneda lemma, the functor $E^\bA : \Stab_R \to \Spt$ is given by
  \begin{equation*}
    E^\bA(\bA') = E(\bA \otimes \bA')
  \end{equation*}
for every $\bA' \in \Stab_R$.
\end{rem}

\begin{defn}
  We say that $E$ is \emph{additive} if it sends split exact sequences of stable \inftyCats \cite[Def.~5.18]{BlumbergGepnerTabuada} to split exact triangles of spectra.
  We say $E$ is \emph{localizing} if it sends short exact sequences of stable \inftyCats \cite[Def.~5.12]{BlumbergGepnerTabuada} to exact triangles of spectra.
  Note that, unlike \cite{BlumbergGepnerTabuada}, we do not require that $E$ commutes with filtered colimits.
\end{defn}

\begin{rem}
Note that if $E$ is additive (resp. localizing), then the same holds for $E^\bA$ for any $\bA \in \Stab_R$.
\end{rem}

\begin{exam}\label{exam:K}
  Let $\K : \Stab \to \Spt$ denote the algebraic K-theory functor.
  Recall that this is defined using the Waldhausen $S_\bullet$-construction (see \cite[Rmk.~1.2.2.5]{HA-20170918}, \cite[Def.~7.1]{BlumbergGepnerTabuada}, or \cite[Sect.~10]{Barwick}) and takes values in connective spectra.
  Then $\K$ is additive by Waldhausen's additivity theorem (\cite[Prop.~7.10]{BlumbergGepnerTabuada}).

  Let $\mathbb{K} : \Stab \to \Spt$ denote the nonconnective algebraic K-theory functor, defined e.g. as in \cite[Def.~9.6]{BlumbergGepnerTabuada} following Schlichting.
  This is a localizing invariant such that $\mathbb{K}_{\geq 0} \simeq \K$.
\end{exam}


\ssec{The projective bundle formula}
\label{ssec:bass/projbun}

In this subsection we prove a projective bundle formula computing $E(\P^{1,\flat}_{R})$ for the \emph{flat} projective line over any \cnEring $R$, and any $R$-linear additive invariant $E$.
This essentially follows from a result of Lurie \cite[Thm.~7.2.2.1]{SAG-20180204}.

\sssec{}

Consider the following subsets of $\bZ \times \bZ$:
\begin{itemize}
  \item $M^+$ consists of pairs $(m,n)$ with $m+n=0$ and $m\ge 0$.
  \item $M^-$ consists of pairs $(m,n)$ with $m+n=0$ and $n\ge 0$.
  \item $M^\pm$ consists of pairs $(m,n)$ with $m+n=0$.
\end{itemize}
We view each of these as (additive) discrete commutative monoids.
For any \cnEring $R$, we write $R[T]$, $R[T^{-1}]$ and $R[T^{\pm}]$ for the monoid algebras $R \otimes \Sigma^\infty_+(M^+)$, $R \otimes \Sigma^\infty_+(M^-)$, and $R \otimes \Sigma^\infty_+(M^\pm)$, respectively.

\sssec{}

We write $p_+ : \Spec(R[T]) \to \Spec(R)$, $p_- : \Spec(R[T^{-1}]) \to \Spec(R)$, and $p_{\pm} : \Spec(R[T^\pm]) \to \Spec(R)$ for the respective projections.
Note that under the obvious isomorphisms $M^+ \simeq \bN \simeq M^-$, both $\Spec(R[T])$ and $\Spec(R[T^{-1}])$ are canonically identified with the \emph{flat affine line} $\A^{1,\flat}_R$ (\defnref{defn:flat affine space}).
Similarly, under the isomorphism $M^\pm \simeq \bZ$, $R[T^{\pm}]$ is identified with the monoid algebra $R \otimes \Sigma^\infty_+(\bZ)$.
It can also be identified with the localization of $R[T]$ away from $T \in \pi_0(R[T]) \simeq \pi_0(R)[T]$, or the localization of $R[T^{-1}]$ away from $T^{-1}$.
In particular, the projections $p_+$, $p_-$ and $p_\pm$ are fibre-smooth in the sense of \cite[Def.~11.2.3.1]{SAG-20180204}, and the affine spectral scheme $\bG^{\flat}_{m,R} = \Spec(R[T^{\pm}])$ is equipped with open immersions
  \begin{equation*}
    j_+ : \Spec(R[T^{\pm}]) \to \Spec(R[T]),
    \quad
    j_- : \Spec(R[T^{-}]) \to \Spec(R[T]).
  \end{equation*}

The construction of the zero and unit sections of $\A^{1,\flat}$ described in \remref{rem:A^n_flat structure} can be adapted as follows.
We recall the constructions of the zero and unit sections of $\Spec(R[T])$.
Consider the set $\{0, 1\}$, viewed as a pointed multiplicative monoid with base point $0$ and identity element $1$.
Since $M^+$ is freely generated as a (discrete) commutative monoid by the element $(1,0) \in M^+$, either choice of element $i \in \{0,1\}$ determines a unique homomorphism $M^+ \to \{0,1\}$ sending $(1,0) \mapsto i$.
Each of these gives rise to $\Einfty$-ring homomorphisms
  \begin{equation*}
    \sigma_i : R[T] \to R\otimes \Sigma^\infty(\{0,1\}) \simeq R,
  \end{equation*}
where we identify $\{0,1\}$ with the pointed $0$-sphere $S^0$.
We let $s_i$ denote the induced morphisms
  \begin{equation*}
    s_i : \Spec(R) \to \Spec(R[T])
  \end{equation*}
for each $i\in\{0,1\}$.
The obvious analogous construction gives sections $s_i : \Spec(R) \to \Spec(R[T^{-1}])$.
Similarly it is clear that the unit section $s_1$ factors through a morphism $s_1 : \Spec(R) \to \Spec(R[T^{\pm}])$.

\sssec{}

We denote by $\P^{1,\flat}_R$ the \emph{flat projective line} over $R$, see \cite[Constr.~5.4.1.3]{SAG-20180204} (where it is denoted $\P^1_R$).
This is equipped with a canonical morphism $q : \P^{1,\flat}_{R} \to \Spec(R)$ which is fibre-smooth.
Moreover there is a cartesian and cocartesian square of spectral schemes
  \begin{equation}\label{eq:P^1 cover}
    \begin{tikzcd}
      \Spec(R[T^\pm]) \ar{r}{j_-}\ar{d}{j_+}
        & \Spec(R[T^{-1}]) \ar{d}{k_-}
      \\
      \Spec(R[T]) \ar{r}{k_+}
        & \P^{1,\flat}_R
    \end{tikzcd}
  \end{equation}
where every arrow is an open immersion.
In particular, this is a Nisnevich square.

\begin{notat}
For a functor $E : \Stab_R \to \Spt$, we write:
  \begin{equation*}
    E^+ := E^{\Perf(R[T])},
    \quad
    E^- := E^{\Perf(R[T^{-1}])},
    \quad
    E^\pm := E^{\Perf(R[T^{\pm}])},
    \quad
    E^\boxtimes := E^{\Perf(\P^{1,\flat}_R)}.
  \end{equation*}
Note that $E^+(\bA) = E(\bA\otimes\Perf(R[T]))$ for $\bA \in \Stab_R$, and similarly for $E^-$, $E^\pm$ and $E^\boxplus$.
\end{notat}

\begin{thm}\label{thm:projbun}
Let $E$ be an additive invariant of $R$-linear stable \inftyCats.
Then the two functors $\Perf(\Spec(R)) \to \Perf(\P^{1,\flat}_R)$, given by $\sF \mapsto q^*(\sF)$ and $\sF \mapsto q^*(\sF) \otimes \sO(-1)$, induce a canonical isomorphism
  \begin{equation*}
    q^* \oplus (q^*\otimes\sO(-1)) : E \oplus E \isoto E^\boxplus.
  \end{equation*}
\end{thm}

\begin{proof}
By Yoneda, it will suffice to show that the maps
  \begin{equation*}
    E(\bA) \oplus E(\bA) \to E(\bA \otimes \Perf(\P^{1,\flat}_{R}))
  \end{equation*}
are invertible for every $\bA \in \Stab_R$.
For this it will suffice to show that the map
  \begin{equation*}
    E(R) \oplus E(R) \to E(\P^{1,\flat}_{R})
  \end{equation*}
is invertible for every additive invariant $E$ (as we can then apply this to every $E^\bA$).

By \cite[Thm.~7.2.2.1]{SAG-20180204} there is a semi-orthogonal decomposition on $\Qcoh(\P^{1,\flat}_{R})$ into two full subcategories both canonically equivalent to $\Qcoh(\Spec(R)) \simeq \Mod_R$.
Moreover, an inspection of the proof of \emph{loc.~cit.} shows that this restricts to a semi-orthogonal decomposition on $\Perf(\P^{1,\flat}_{R})$ by two full subcategories both equivalent to $\Perf(\Spec(R)) \simeq \Perf_R$.
Indeed, both functors $q^*$ and $q_*$ preserve perfect complexes (the latter by \cite[Thm.~6.1.3.2]{SAG-20180204}).
The claim then follows by definition of additive invariants.
\end{proof}

\begin{rem}\label{rem:alpha projbun}
Let $\alpha$ denote the composite morphism
  \begin{equation*}
    \alpha : E \oplus E
      \xrightarrow{\mu} E \oplus E
      \xrightarrow{q^* \oplus (q^*\otimes\sO(-1))} E^\boxplus,
  \end{equation*}
where $\mu$ is the isomorphism induced by the invertible matrix
  \begin{equation*}
    \begin{bmatrix}
      1 & 1\\
      0 & -1
    \end{bmatrix}.
  \end{equation*}
By construction, $\alpha$ fits in the commutative diagram
  \begin{equation}\label{eq:alpha with k^*}
    \begin{tikzcd}
        & E\oplus E\ar{d}{\alpha}\ar{rd}{p_+^*\oplus 0}\ar[swap]{ld}{p_-^*\oplus 0}
        &
      \\
      E^-
        & E^\boxtimes\ar[swap]{r}{k_+^*}\ar{l}{k_-^*}
        & E^+.
    \end{tikzcd}
  \end{equation}
\end{rem}


\ssec{The Bass fundamental sequence}
\label{ssec:bass/bassseq}

\begin{thm}\label{thm:bassseq}
Let $E : \Stab_R \to \Spt$ be a localizing invariant.
Then for any $R$-linear stable \inftyCat $\bA$, there is a natural exact sequence of abelian groups
  \begin{equation*}
    0 \to E_n^\bA(R)
      \xrightarrow{(p_+^*,-p_-^*)} E_n^\bA(R[T]) \oplus E_n^\bA(R[T^{-1}])
      \xrightarrow{(j_+^*,j_-^*)} E_n^\bA(R[T^\pm])
      \xrightarrow{\partial} E_{n-1}^\bA(R)
      \to 0
  \end{equation*}
for every integer $n$, where we write $E_n = \pi_n E$.
\end{thm}

\begin{proof}
By replacing $E$ with $E^\bA$, we may assume that $\bA = \Perf_R$.
Since $E$ satisfies Nisnevich descent (e.g. \cite[App.~A]{ClausenMathewNaumannNoelDescent}), the Nisnevich square \eqref{eq:P^1 cover} gives rise to a cartesian square
  \begin{equation*}
    \begin{tikzcd}
      E(\P^{1,\flat}_R) \ar{r}{k_-^*}\ar{d}{k_+^*}
        & E(R[T^{-1}]) \ar{d}{j_-^*}
      \\
      E(R[T]) \ar{r}{j_+^*}
        & E(R[T, T^{-1}])
    \end{tikzcd}
  \end{equation*}
and hence to a Mayer--Vietoris long exact sequence
  \begin{equation}\label{eq:P^1 MV}
    \begin{split}
      \cdots
        \to E_{n+1}(&R[T,T^{-1}])
        \xrightarrow{\partial} E_n(\P^{1,\flat}_R)
        \xrightarrow{(k_+^*, -k_-^*)}\\ &E_n(R[T]) \oplus E_n(R[T^{-1}])
        \xrightarrow{j_+^* \oplus j_-^*} E_n(R[T,T^{-1}])
        \xrightarrow{\partial} \cdots
    \end{split}
  \end{equation}
By the projective bundle formula (\thmref{thm:projbun}) and \remref{rem:alpha projbun}, there is a canonical isomorphism $E(R) \oplus E(R) \simeq E(\P^{1,\flat}_{R})$ under which the map $(k_+^*, -k_-^*)$ in \eqref{eq:P^1 MV} is $(p_+^*, p_-^*)$ on the first copy of $E_n(R)$ and $(0,0)$ on the second.
Here $p_+$ and $p_-$ denote the respective projections $p_+ : \Spec(R[T]) \to \Spec(R)$ and $p_- : \Spec(R[T^{-1}]) \to \Spec(R)$.
Since the boundary map is then surjective onto the second copy of $E_n(R)$, and the zero section induces canonical retractions of both $p_+^*$ and $p_-^*$, we see that the long exact sequence \eqref{eq:P^1 MV} splits up into short exact sequences as in the claim.
\end{proof}

\begin{rem}\label{rem:Bott class}
In the case where the localizing invariant $E$ is nonconnective algebraic K-theory $\mathbb{K}$ (\examref{exam:K}), the map $\partial : \mathbb{K}_n(R[T^\pm]) \to \mathbb{K}_{n-1}(R)$ in the Bass fundamental sequence admits a natural \emph{splitting}, up to an automorphism of $\mathbb{K}_{n-1}(R)$.
Indeed, consider the automorphism of $R[T,T^{-1}]$ given by multiplication by $T$.
This induces a point $b \in \mathbb{K}(R[T,T^{-1}])[-1]$ which we call the \emph{Bott class}.
Now cup product with $b$ induces a canonical map
  \begin{equation*}
    \mathbb{K}(R)
      \xrightarrow{p_\pm^*} \mathbb{K}(R[T^\pm])
      \xrightarrow{b\cup} \mathbb{K}(R[T^\pm])[-1]
      \xrightarrow{\partial} \mathbb{K}(R)
  \end{equation*}
which we claim is invertible.
By $\mathbb{K}(R)$-linearity, this is equivalent to the assertion that $\partial$ sends $b\in \mathbb{K}_1(R[T^\pm]) \simeq \K_1(R[T^\pm])$ to a unit in $\mathbb{K}_0(R) \simeq \K_0(R)$.
Since the $1$-truncation $\tau_{\le 1}(\K)$ is insensitive to positive homotopy groups \cite[Lect.~20, Cor.~4]{LurieKLectures}, we may replace $R$ by $\pi_0(R)$.
Now the claim is classical, see e.g. the proof of \cite[Thm.~6.1(b)]{ThomasonTrobaugh}.
\end{rem}


\ssec{Delooping localizing invariants of connective spectra}

\begin{notat}
We let $\Spt_{\ge 0}$ denote the full subcategory of $\Spt$ spanned by connective spectra.
A \emph{connective fibre sequence} of spectra is a diagram
  \begin{equation*}
    F \to X \to Y,
  \end{equation*}
together with a null-homotopy of the composite $F \to Y$, such that the induced map $F \to \Fib(X \to Y)$ induces an isomorphism
  \begin{equation*}
    \tau_{\ge 0}(F) \simeq \tau_{\ge 0}(\Fib(X \to Y))
  \end{equation*}
of connective spectra.
\end{notat}

\begin{defn}\label{defn:localizing connective}
  We say that a functor $\Stab_R \to \Spt_{\ge 0}$ is \emph{localizing} if it sends short exact sequences to connective fibre sequences.
\end{defn}

\begin{thm}\label{thm:deloop localizing invariant}
Let $R$ be a connective $\Einfty$-ring.
The assignment $E \mapsto \tau_{\ge 0}(E)$ determines an equivalence
  \begin{equation*}
    \Fun_{\mrm{loc}}(\Stab_R, \Spt) \to \Fun_{\mrm{loc}}(\Stab_R, \Spt_{\ge 0})
  \end{equation*}
from the \inftyCat of localizing invariants $\Stab_R \to \Spt$, to the \inftyCat of localizing invariants $\Stab_R \to \Spt_{\ge 0}$.
\end{thm}

\begin{exam}\label{exam:K^B}
  Let $\K : \Stab \to \Spt_{\ge 0}$ and $\mathbb{K} : \Stab \to \Spt$ be as in \examref{exam:K}.
  Since $\mathbb{K}$ is localizing and satisfies $\mathbb{K}_{\geq 0} \simeq \K$, it follows from \thmref{thm:deloop localizing invariant} that there is a canonical isomorphism of localizing invariants
  \[
    \mathbb{K} \simeq \K^\B.
  \]
\end{exam}


\ssec{The Bass construction}
\label{ssec:bass/bassconstr}

Let $E : \Stab_R \to \Spt$ be an arbitrary functor.

\sssec{}
Define $V(E)$ and $W(E)$ such that there are cocartesian squares
  \begin{equation}\label{eq:V and W diagram}
    \begin{tikzcd}
      E\ar{r}{\mrm{incl}_1}\ar{d}
        & E\oplus E\ar{r}{\alpha}\ar{d}{\mrm{pr}_2}
        & E^\boxplus\ar{r}{(k_+^*,-k_-^*)}\ar{d}
        & E^+\oplus E^-\ar{d}
      \\
      0\ar{r}
        & E\ar{r}
        & W(E)\ar{r}
        & V(E)
    \end{tikzcd}
  \end{equation}
in $\Fun(\Stab_R,\Spt)$.
Here $\alpha$ is as in \remref{rem:alpha projbun}.
We denote by $\psi_E$ the composite $E \to W(E) \to V(E)$.

\begin{rem}\label{rem:psi=0}
The commutative diagram \eqref{eq:alpha with k^*} provides a null-homotopy of the composite
  \begin{equation*}
    E \xrightarrow{\mrm{incl}_2} E\oplus E
      \xrightarrow{\alpha} E^\boxplus
      \xrightarrow{(k_+^*,-k_-^*)} E^+\oplus E^-.
  \end{equation*}
Since the composite $\mrm{pr}_2\circ \mrm{incl}_2$ is the identity, combining this with the commutative diagram \eqref{eq:V and W diagram} yields a canonical null-homotopy of the morphism $\psi_E : E \to V(E)$.
\end{rem}

\begin{rem}
The commutative diagram \eqref{eq:alpha with k^*} identifies the upper horizontal composite with the morphism $(p_+^*,-p_-^*) : E \to E^+\oplus E^-$.
It follows that $V(E)$ fits in an exact triangle
  \begin{equation*}
    E \xrightarrow{(p_+^*, -p_-^*)} E^+ \oplus E^- \to V(E).
  \end{equation*}
Moreover, since the morphisms $p_+^*$ and $p_-^*$ admit splittings induced by the homomorphisms $R[T] \to R$ and $R[T^{-1}] \to R$, respectively, it follows that the associated long exact sequence splits into short exact sequences
  \begin{equation}\label{eq:V(E) ses}
    0 \to \pi_n(E) \xrightarrow{(p_+^*, -p_-^*)} \pi_n(E^+) \oplus \pi_n(E^-) \to \pi_n(V(E)) \to 0
  \end{equation}
for every $n$.
\end{rem}

\begin{rem}
Note that we have commutative squares
  \begin{equation*}
    \begin{tikzcd}
      E^\boxplus \ar{r}{k_+^*}\ar{d}{k_-^*}
        & E^+ \ar{d}{j_+^*}
      \\
      E^- \ar{r}{j_-^*}
        & E^\pm,
    \end{tikzcd}
    \qquad
    \begin{tikzcd}
      E^\boxplus \ar{r}{(k_+^*,-k_-^*)}\ar{d}
        & E^+\oplus E^- \ar{d}{j_+^*\oplus j_-^*}
      \\
      0 \ar{r}
        & E^\pm,
    \end{tikzcd}
  \end{equation*}
where the right-hand square is induced from the left-hand one.
This gives rise to a canonical morphism $\theta_E : V(E) \to E^\pm$ fitting into the commutative diagram:
  \begin{equation*}
    \begin{tikzcd}
        E^\boxplus\ar{r}{(k_+^*,-k_-^*)}\ar{d}
        & E^+\oplus E^-\ar{d}\ar[bend left=60, looseness=1.2]{dd}{j_+^*\oplus j_-^*}
      \\
        W(E)\ar{r}\ar{d}
        & V(E)\ar[dashed]{d}{\theta_E}
      \\
        0\ar{r}
        & E^\pm
    \end{tikzcd}
  \end{equation*}
\end{rem}

\begin{constr}\label{constr:Bass}
Let $E : \Stab_R \to \Spt$ be a functor.
Denote by $U(E)$ the fibre of the morphism $\theta_E : V(E) \to E^\pm$, so that there is a fibre sequence
  \begin{equation*}
    U(E) \to V(E) \xrightarrow{\theta_E} E^\pm.
  \end{equation*}
The canonical null-homotopy of the composite $W(E) \to V(E) \xrightarrow{\theta_E} E^\pm$ defined above gives rise to a canonical morphism $W(E) \to U(E)$.
In particular, we get a canonical morphism
  \begin{equation*}
    \phi_E : E \to W(E) \to U(E).
  \end{equation*}
This gives rise to a tower
  \begin{equation*}
    E \xrightarrow{\phi_E} U(E) \xrightarrow{U(\phi_{E})} U^2(E) \xrightarrow{U^2(\phi_{E})} \cdots
  \end{equation*}
whose colimit we denote $E^\B$, and call the \emph{Bass construction} on $E$.
\end{constr}

\begin{rem}\label{rem:U(phi)=phi_U(E)}
Since the functors $V$ and $(-)^\pm$ commute with colimits and with $(-)^\bA$ for $\bA\in\Stab_R$, the same holds for $U$.
In particular, it follows that there are canonical identifications $U^k(\phi_E) \simeq \phi_{U^k(E)}$ for each $k\ge 0$.
\end{rem}

\begin{thm}\label{thm:Bass}
Let $E : \Stab_R \to \Spt_{\ge 0}$ be a functor.
The Bass construction $E^\B$ satisfies the following properties:
\begin{thmlist}
  \item\label{thm:Bass/deloops}
If $E$ satisfies the projective bundle formula and is $Q_\pm$-excisive, then the natural morphism $E \to E^\B$ induces an isomorphism $E \simeq \tau_{\ge0}(E^\B)$.

  \item\label{thm:Bass/localizing}
If $E : \Stab_R \to \Spt_{\ge 0}$ is a localizing invariant, then the functor $E^\B : \Stab_R \to \Spt$ is a localizing invariant.
\end{thmlist}
\end{thm}


\ssec{Proof of \thmref{thm:Bass}}

\begin{lem}\label{lem:bass hilfsatz}
Let $E : \Stab_R \to \Spt_{\ge0}$ be a functor.
If $E$ satisfies the projective bundle formula and is $Q_\pm$-excisive, then we have:
\begin{thmlist}
  \item\label{item:bass hilfsatz/1}
There is a connective fibre sequence
  \begin{equation*}
    E \xrightarrow{\psi_E} V(E) \xrightarrow{\theta_E} E^\pm
  \end{equation*}
which is natural in $E$.
  \item\label{item:bass hilfsatz/2}
There exist canonical morphisms
  \begin{equation*}
    \tau_{\ge0}(\Omega(E^\pm)) \to E,
    \quad
    \sigma_E : E \to \tau_{\ge0}(\Omega(E^\pm))
  \end{equation*}
which exhibit $E$ as a retract of $\tau_{\ge0}(\Omega(E^\pm))$, and are natural in $E$.
\end{thmlist}
\end{lem}

\begin{proof}
We first show the following weaker version of \itemref{item:bass hilfsatz/1}:
\begin{enumerate}[label={($\ast$)}]
  \item\label{item:bass hilfsatz/1bis}
There is a connective fibre sequence
  \begin{equation*}
    \tau_{\ge0}(\Omega(E^\pm)) \to E \xrightarrow{\psi_E} V(E).
  \end{equation*}
\end{enumerate}
Since $E$ satisfies the projective bundle formula, $\alpha : E\oplus E\to E^\boxplus$ is invertible.
Considering the diagram \eqref{eq:V and W diagram}, it follows that $E\to W(E)$ is also invertible and that it will suffice to show that the diagram
  \begin{equation*}
    \tau_{\ge0}(\Omega(E^\pm)) \to W(E) \to V(E)
  \end{equation*}
is a connective fibre sequence.
Since the \inftyCat $\Spt_{\ge0}$ is prestable, the right-hand square in \eqref{eq:V and W diagram} is also cartesian \cite[Cor.~C.1.2.6]{SAG-20180204}, and in particular induces an isomorphism on fibres.
The fact that $E$ is $Q_\pm$-excisive implies that the fibre\footnotemark~of the upper arrow $(k_+^*,-k_-^*) : E^\boxplus \to E^+\oplus E^-$ is $\tau_{\ge0}(\Omega(E^\pm))$.
\footnotetext{Note that this fibre, computed in $\Spt_{\ge0}$, is the same as the connective cover of the fibre computed in $\Spt$.}
This shows claim~\itemref{item:bass hilfsatz/1bis}.

For part~\itemref{item:bass hilfsatz/2}, the canonical null-homotopy of the morphism $\psi_E : E \to V(E)$ (\remref{rem:psi=0}) gives rise to a morphism $\sigma_E$ fitting in the commutative diagram
  \begin{equation*}
    \begin{tikzcd}
      E \ar[dashed,swap]{d}{\sigma_E}\ar[equals]{rd}\ar[bend left]{rrd}
        & 
        &
      \\
      \tau_{\ge0}(\Omega(E^\pm)) \ar{r}
        & E \ar{r}{\psi_E}
        & V(E),
    \end{tikzcd}
  \end{equation*}
since the lower row is a fibre sequence by claim~\itemref{item:bass hilfsatz/1bis}.

Finally we prove \itemref{item:bass hilfsatz/1}.
By \itemref{item:bass hilfsatz/2}, the diagram $E \to V(E) \to E^\pm$ is a retract of the diagram
  \begin{equation*}
    \tau_{\ge0}(\Omega(E^\pm))
      \to \tau_{\ge0}(\Omega(V(E)^\pm))
      \to \tau_{\ge0}(\Omega((E^\pm)^\pm)).
  \end{equation*}
It will suffice to show that the latter is a fibre sequence.
Since the functor $E \mapsto E^\pm$ is left-exact and commutes with $\tau_{\ge0}$, this follows from the fact that the diagram
  \begin{equation*}
    \tau_{\ge0}(\Omega(E))
      \to \tau_{\ge0}(\Omega(V(E)))
      \to \tau_{\ge0}(\Omega(E^\pm))
  \end{equation*}
is a fibre sequence, by \itemref{item:bass hilfsatz/1bis}.
\end{proof}

\begin{cor}\label{cor:U deloops}
Let $E : \Stab_R \to \Spt_{\ge 0}$ be a localizing invariant.
Then the morphism $\phi_E : E \to U(E)$ induces an isomorphism $E \simeq \tau_{\ge0}(U(E))$.
\end{cor}

\begin{proof}
The claim is equivalent to the assertion that the diagram
  \begin{equation*}
    E \to V(E) \to E^\pm
  \end{equation*}
is a connective fibre sequence.
This is \lemref{lem:bass hilfsatz}\itemref{item:bass hilfsatz/1}, which applies since $E$ is localizing.
\end{proof}

\begin{lem}\label{lem:bass hilfsatz nonconnective}
Let $E : \Stab_R \to \Spt$ be a functor.
If $E$ is localizing, then we have:
\begin{thmlist}
  \item\label{item:bass hilfsatz nonconnective/1}
There is an exact triangle
  \begin{equation*}
    E \xrightarrow{\psi_E} V(E) \xrightarrow{\theta_E} E^\pm
  \end{equation*}
which is natural in $E$.
In other words, the morphism $\phi_E : E \to U(E)$ is invertible.
  \item\label{item:bass hilfsatz nonconnective/2}
The canonical morphism $E \to E^\B$ is invertible.
  \item\label{item:bass hilfsatz nonconnective/3}
There exist canonical morphisms
  \begin{equation*}
    \Omega(E^\pm) \to E,
    \quad
    \sigma_E : E \to \Omega(E^\pm)
  \end{equation*}
which exhibit $E$ as a retract of $\Omega(E^\pm)$, and are natural in $E$.
\end{thmlist}
\end{lem}

\begin{proof}
Parts (i) and (iii) follow by the same argument as in the proof of \lemref{lem:bass hilfsatz}.
The second follows from (i).
\end{proof}

\begin{lem}\label{lem:U preserves tau>=0-equivalences}
The functor $E \mapsto U(E)$ preserves $\tau_{\ge0}$-equivalences.
That is, let $E \to E'$ be a morphism in $\Fun(\Stab_R,\Spt)$ which induces an isomorphism $\tau_{\ge 0}(E) \simeq \tau_{\ge0}(E')$.
Then the induced map
  \begin{equation*}
    \tau_{\ge 0}(U(E)) \to \tau_{\ge 0}(U(E'))
  \end{equation*}
is invertible.
\end{lem}

\begin{proof}
Note that the analogous property holds for the functors $(-)^\pm$ and $V(-)$: in fact, they are even t-exact (the latter in view of the exact sequences \eqref{eq:V(E) ses}).
As $U(-)$ is the fibre of the morphism $V(-) \to (-)^\pm$, the claim follows.
\end{proof}

\begin{lem}\label{lem:bass hilfsatz for E^B}
Let $E : \Stab_R \to \Spt$ be a functor.
We have:
\begin{thmlist}
  \item\label{item:bass hilfsatz for E^B/1}
There is an exact triangle
  \begin{equation*}
    E^\B \xrightarrow{\psi_{E^\B}} V(E^\B) \xrightarrow{\theta_{E^\B}} (E^\B)^\pm
  \end{equation*}
which is natural in $E$.
  \item\label{item:bass hilfsatz for E^B/2}
There exist canonical morphisms
  \begin{equation*}
    \Omega((E^\B)^\pm) \to E^\B,
    \quad
    \sigma_{E^\B} : E^\B \to \Omega((E^\B)^\pm)
  \end{equation*}
which exhibit $E^\B$ as a retract of $\Omega((E^\B)^\pm)$, and are natural in $E$.
\end{thmlist}
\end{lem}

\begin{proof}
Part \itemref{item:bass hilfsatz for E^B/2} will follow from \itemref{item:bass hilfsatz for E^B/1} as in the proof of \lemref{lem:bass hilfsatz}.
For \itemref{item:bass hilfsatz for E^B/1} it will suffice to show that the morphism $\phi_{E^\B} : E^\B \to U(E^\B)$ is invertible.
Since the functors $U$, $V$ and $W$ each commute with colimits, it is clear from the construction that this morphism is the colimit of the morphisms $\phi_{U^k(E)} \simeq U^k(\phi_E) : U^k(E) \to U^{k+1}(E)$, $k\ge 0$ (\remref{rem:U(phi)=phi_U(E)}).
That is, $\phi_{E^\B}$ is identified with the canonical morphism fitting in the diagram
  \begin{equation*}
    \begin{tikzcd}
      E \ar{r}{\phi_E}\ar{d}{\phi_E}
        & U(E) \ar{r}{U(\phi_E)}\ar{d}{U(\phi_E)}
        & U^2(E) \ar{r}{U^2(\phi_E)}\ar{d}{U^2(\phi_E)}
        & \cdots \ar{r}
        & E^\B \ar{d}{\phi_{E^\B}}
      \\
      U(E) \ar[swap]{r}{U(\phi_E)}
        & U^2(E) \ar[swap]{r}{U^2(\phi_E)}
        & U^3(E) \ar[swap]{r}{U^3(\phi_E)}
        & \cdots \ar[swap]{r}
        & U(E^\B),
    \end{tikzcd}
  \end{equation*}
and is clearly invertible.
\end{proof}

\begin{cor}\label{cor:E^B retract}
Let $E : \Stab_R \to \Spt$ be a functor.
For every integer $n\ge 0$, denote by $(E^\B)^{\pm n} : \Stab_R \to \Spt$ the functor
  \begin{equation*}
    ((((E^\B)^{\pm})^\pm)^\cdots)^\pm
  \end{equation*}
obtained from the Bass construction $E^\B$ by an $n$-fold iteration of the functor $(-)^\pm$ (e.g. $(E^\B)^{\pm 0} = E^\B$).
Then for each $n\ge 0$, the functor
  \begin{equation*}
    \Sigma^{n}(E^\B) : \Stab_R \to \Spt
  \end{equation*}
is a retract of $(E^\B)^{\pm n} : \Stab_R \to \Spt$.
\end{cor}

\begin{proof}
By induction, it will suffice to consider $n=1$ and, by adjunction, to show that $E^\B$ is a retract of $\Omega(E^\B)^\pm$.
This is \lemref{lem:bass hilfsatz for E^B}\itemref{item:bass hilfsatz for E^B/2}.
\end{proof}

\sssec{Proof of \thmref{thm:Bass}\itemref{thm:Bass/deloops}}
By \corref{cor:U deloops}, $\phi_E : E \to U(E)$ is a $\tau_{\ge0}$-equivalence.
It follows then from \lemref{lem:U preserves tau>=0-equivalences} that $U^k(\phi_E) : U^k(E) \to U^{k+1}(E)$ is a $\tau_{\ge0}$-equivalence for every $k\ge 0$.
It follows that the transfinite composite $E \to E^\B$ is also a $\tau_{\ge0}$-equivalence.

\sssec{Proof of \thmref{thm:Bass}\itemref{thm:Bass/localizing}}

Let $\bA' \to \bA \to \bA''$ be an exact sequence of small stable \inftyCats and consider the induced diagram of spectra
  \begin{equation*}
    E^\B(\bA') \to E^\B(\bA) \to E^\B(\A'').
  \end{equation*}
To show that this is an exact triangle, it will suffice to show that each of the induced diagrams
  \begin{equation*}
    \tau_{\ge n}(E^\B(\bA')) \to \tau_{\ge n}(E^\B(\bA)) \to \tau_{\ge n}(E^\B(\A''))
  \end{equation*}
is a fibre sequence in $\Spt_{\ge n}$, for every $n\le 0$.
In other words, it will suffice to show that each functor $\tau_{\ge n}(E^\B) : \Stab_R \to \Spt_{\ge n}$ is a localizing invariant.

By \corref{cor:E^B retract} we know that $\tau_{\ge n}(E^\B)$ is a retract of $\tau_{\ge 0}((E^\B)^{\pm n})$.
Since $(-)^\pm$ is left-exact, the latter is isomorphic to $\tau_{\ge0}(E^\B)^{\pm n} \simeq E^{\pm n}$ (\thmref{thm:Bass}\itemref{thm:Bass/deloops}), which is localizing because $E$ is.


\ssec{Proof of \thmref{thm:deloop localizing invariant}}

We are now ready to prove \thmref{thm:deloop localizing invariant}, which asserts that the canonical functor
  \begin{equation*}
    \tau_{\ge0} : \Fun_{\mrm{loc}}(\Stab_R, \Spt) \to \Fun_{\mrm{loc}}(\Stab_R, \Spt_{\ge 0}),
  \end{equation*}
given by the assignment $E \mapsto \tau_{\ge 0}(E)$, is an equivalence.

\sssec{}
Note that the Bass construction (\constrref{constr:Bass}) defines a functor $E \mapsto E^\B$ from $\Fun_{\mrm{loc}}(\Stab_R, \Spt_{\ge0})$ to $\Fun_{\mrm{loc}}(\Stab_R, \Spt)$.
We claim first that this is a fully faithful left adjoint to $\tau_{\ge 0}$.
For $E \in \Fun_{\mrm{loc}}(\Stab_R, \Spt_{\ge0})$, the unit map
  \begin{equation*}
    \eta_E : E \isoto \tau_{\ge 0}(E^\B)
  \end{equation*}
comes from the natural isomorphisms of \thmref{thm:Bass}\itemref{thm:Bass/deloops}.
For $E \in \Fun_{\mrm{loc}}(\Stab_R, \Spt)$, the co-unit map
  \begin{equation*}
    \varepsilon_E : (\tau_{\ge 0}(E))^\B \to E^\B
  \end{equation*}
is induced from $\tau_{\ge0}(E) \to E$ in view of the fact that the natural map $E \to E^\B$ is invertible (\lemref{lem:bass hilfsatz nonconnective}\itemref{item:bass hilfsatz nonconnective/2}).
One easily verifies the triangle identities.

\sssec{}

In order to conclude that the functor $\tau_{\ge 0}$ is an equivalence, it will suffice to show that it is conservative (so that the co-unit maps are necessarily isomorphisms).
Let $E \to E'$ be a morphism of localizing invariants in $\Fun_{\mrm{loc}}(\Stab_R, \Spt)$, and suppose that the induced morphism $\tau_{\ge0}(E) \to \tau_{\ge0}(E')$ is invertible.
We argue by decreasing induction on $n$ that the map $\pi_n(E(\bA)) \to \pi_n(E'(\bA))$ is invertible for every $n\le 0$ and every $\bA \in \Stab_R$.
For $n=0$ this holds by assumption.
The induction step follows from the Bass fundamental sequence (\thmref{thm:bassseq}).
The claim follows.


\section{Homotopy invariant K-theory}
\label{sec:kh}


\ssec{Homotopy invariant K-theory}
\label{ssec:kh/kh}

\begin{notat}
  Given a spectral scheme $S$, let $\Perf(S)$ denote the stable \inftyCat of perfect complexes on $S$.
  Denote by $\K(S)$ and $\KB(S)$ the spectra
  \[
    \K(\Perf(S)), \quad
    \KB(\Perf(S)),
  \]
  respectively.
  Here $\K$ is algebraic K-theory (\examref{exam:K}) and $\KB$ is the Bass construction, equivalent to nonconnective algebraic K-theory $\mathbb{K}$ (\examref{exam:K^B}).
\end{notat}

\begin{constr}\label{constr:KH}
  For any spectral scheme $S$, consider the spectrum
  \[
    \KH(S) := \colim \K(S \times \A^n),
  \]
  where $\A^n$ is the $n$-dimensional spectral affine space (\examref{exam:A^n}) and the colimit is indexed by the opposite of the (cosifted) full subcategory $\bA_S \sub \Aff_{/S}$ whose objects are spectral affine spaces $S \times \A^n$ ($n\ge 0$).
  We also write $\KH(R) = \KH(\Spec(R))$ for any \cnEring $R$.
\end{constr}

This definition may appear ad-hoc.
In the language of $\sC$-fibred $S^1$-spectra (\examref{exam:motivic S^1-spectra}), we can give a more systematic definition:

\begin{constr}\label{constr:KH fibred}
  Let $\sC_{/S} \sub \Aff_{/S}$ be an admissible subcategory.
  The assignments $X \mapsto \K(X)$ and $X \mapsto \KB(X)$ define presheaves of spectra
  $$(\sC_{/S})^\op \to \Spt,$$
  which we shall denote by $\K|_{\sC}$ and $\KB|_{\sC}$ and which we view as $\sC$-fibred $S^1$-spectra over $S$.
  Then $\KH|_{\sC}$ is the $\sC$-fibred $S^1$-spectrum defined as the $\A^1$-localization of $\KB|_{\sC}$:
  \[ \KH|_{\sC} := \Lhtp(\KB|_{\sC}). \]
  The formula \eqref{eq:formula for Lhtp} shows that the spectrum of global sections recovers $\KH(S)$ as defined above:
  \[ \Gamma(S, \KH|_{\sC}) \simeq \KH(S). \]
\end{constr}

\begin{prop}\label{prop:KH motivic}
  The $\sC$-fibred $S^1$-spectrum $\KH|_{\sC}$ satisfies Nisnevich excision and $\A^1$-homotopy invariance; that is, it is motivic.
  Moreover, the canonical morphism of motivic $\sC$-fibred $S^1$-spectra \[ \L(\KB|_{\sC}) \to \KH|_{\sC} \] is invertible.
\end{prop}

\begin{proof}
Recall that any localizing invariant satisfies Nisnevich excision, see e.g. \cite[Prop.~A.13]{ClausenMathewNaumannNoelDescent}, so $\KB|_{\sC}$ satisfies Nisnevich excision.
Thus by \propref{prop:Lmot of spectra} we have
  $$\L(\KB|_{\sC}) \simeq \Lhtp \LNis(\KB|_{\sC}) \simeq \Lhtp(\KB|_{\sC}),$$
where the latter is $\KH|_{\sC}$ by definition.
\end{proof}

\begin{defn}
  We may repeat \constrref{constr:KH fibred} in the setting of $\sC^\cl$-fibred spectra over $S_\cl$ (notation as in \examref{exam:A^cl}).
  The restriction of $\KB|_{\sC}$ along $u : \sC^\cl_{/S_\cl} \to \sC_{/S}$ (\remref{rem:u}) is the $\sC^\cl$-fibred $S^1$-spectrum $u^*(\KB|_{\sC}) \simeq v_!(\KB|_{\sC})$ which we denote simply by $\KB|_{\sC^\cl}$.
  We define the $\sC^\cl$-fibred $S^1$-spectrum $\KH^\cl|_{\sC^\cl}$ as the $\A^1_\cl$-localization of $\KB|_{\sC^\cl}$:
  \[ \KH^\cl|_{\sC^\cl} := \Lhtpcl(\KB|_{\sC^\cl}). \]
  We define $\KH^\cl(S_\cl)$ as the spectrum of global sections $\Gamma(S_\cl, \KH^\cl|_{\sC^\cl})$.
  This is nothing else than Weibel's homotopy invariant K-theory spectrum (see \cite{cisinski2013descente} for this point of view).
\end{defn}

To formulate the main result of this subsection, we introduce the following notation:

\begin{notat}\label{notat:KHlafds}
  Let $\sA_{/S} \sub \Aff_{/S}$ be a narrow subcategory (e.g. $\sA_{/S} = \Sm_{/S}$), let $\sA^\cl_{/S} \sub \AffCl_{/S}$ be as in \examref{exam:A^cl}, and let $w : \sA_{/S} \to \sA^\cl_{/S_\cl}$ be the classical truncation functor (\constrref{constr:v}).
  Recall the equivalence
  \[
    \L w_! : \MotSpc(\sA_{/S})_\Spt \to \MotSpc(\sA^\cl_{/S_\cl})_\Spt,
    \quad w^* : \MotSpc(\sA^\cl_{/S_\cl})_\Spt \to \MotSpc(\sA_{/S})_\Spt
  \]
  from \thmref{thm:nil-localization V-linear}.
\end{notat}

Then we have:

\begin{thm}\label{thm:v_!KH=KH}
  Let the notation be as in \ref{notat:KHlafds}.
  Then there are canonical isomorphisms
  \[ \L w_!(\KH|_{\sA}) \simeq \KH^\cl|_{\sA^\cl}, \]
  \[ \KH|_{\sA} \simeq w^*(\KH^\cl|_{\sA^\cl}) \]
  of motivic $\sA^\cl$-fibred $S^1$-spectra over $S_\cl$, resp. of motivic $\sA$-fibred $S^1$-spectra over $S$.
\end{thm}

See \ssecref{ssec:kh/coda} for the proof.
Note that this immediately implies \thmref{thm:intro/KH}.

\begin{cor}\label{cor:KH global sections}
  For every quasi-compact quasi-separated spectral algebraic space $S$, there is a canonical isomorphism of spectra $\KH(S) \simeq \KH^\cl(S_\cl)$, functorial in $S$.
\end{cor}

\begin{proof}
  By Nisnevich descent we may assume that $S$ is affine.
  Passing to global sections in \thmref{thm:v_!KH=KH}, we get an isomorphism of spectra
    \begin{equation*}
      \KH(S) = \Gamma(S, \KH|_{\sA})
        \isotoo \Gamma(S, w^*(\KH^\cl|_{\sA^\cl}))
        \simeq \Gamma(S_\cl, \KH^\cl|_{\sA^\cl})
        = \KH^\cl(S_\cl)
    \end{equation*}
  as claimed.
\end{proof}

\begin{cor}
  For any \cnEring $R$, there is a canonical isomorphism of spectra $\KH(R) \simeq \KH^\cl(\pi_0(R))$, functorial in $R$.
\end{cor}


\ssec{Connective comparison}
\label{ssec:kh/unstable}

In this subsection our goal is to prove the following two statements:

\begin{prop}\label{prop:unstable comparison}
  Let the notation be as in \ref{notat:KHlafds}.
  There is a canonical isomorphism of $\sA^\cl$-fibred $S^1$-spectra
  \[
    \L w_! (\K|_{\sA}) \to \L (\K|_{\sA^\cl}).
  \]
\end{prop}

\begin{prop}\label{prop:unstable K Sm-fibred}
  Let the notation be as in \ref{notat:KHlafds}.
  Let $\sB_{/S} \sub \Aff_{/S}$ be a broad subcategory containing $\sA_{/S}$ and let $\iota : \sA_{/S} \hook \sB_{/S}$ denote the inclusion.
  Then the canonical morphisms of $\sB$-fibred $S^1$-spectra
  \begin{align*}
    \LNis \iota_! (\K|_{\sA}) \to \K|_{\sB}\\
    \L \iota_! (\K|_{\sA}) \to \L (\K|_{\sB})
  \end{align*}
  are invertible.
\end{prop}

\begin{cor}\label{cor:LK|B nil-local}
  The motivic $\sB$-fibred $S^1$-spectrum $\L(\K|_{\sB})$ is nil-local (\defnref{defn:nil-local}).
\end{cor}
\begin{proof}
  By \propref{prop:unstable K Sm-fibred}, $\L(\K|_{\sB})$ belongs to the essential image of $\L \iota_! : \MotSpc(\sA_{/S})_\Spt \to \MotSpc(\sB_{/S})_\Spt$.
  Hence it is nil-local by \thmref{thm:nil descent}.
\end{proof}


We will deduce Propositions~\ref{prop:unstable comparison} and \ref{prop:unstable K Sm-fibred} from a representability statement, \propref{prop:X^gp = Omega^infty(K)} below.

\begin{constr}\label{constr:X}
Let $R$ be a \cnEring.
Denote by $\Mod_R^\proj$ the \inftyCat of finitely generated projective $R$-modules, and by $\Mod_R^\free$ the full subcategory of free $R$-modules of finite rank.
Let $X(R)$ denote the underlying \inftyGrpd $(\Mod_R^\free)^\simeq$, obtained by discarding non-invertible morphisms.
Note that $X(R)$ is nothing else than the coproduct of the classifying spaces $\BGL_n(R)$ over $n\ge 0$, where $\GL_n(R)$ is the space of automorphisms of the free $R$-module $R^{\oplus n}$.
Note also that the symmetric monoidal structure on $\Mod_R^\proj$ induces a structure of grouplike $\Einfty$-monoid on $X(R)$.
Moreover, formation of $X(R)$ is functorial and we may regard the assignment $\Spec(R) \mapsto X(R)$ as a presheaf of grouplike $\Einfty$-spaces on the site of affine spectral schemes.
\end{constr}

\begin{prop}\label{prop:X^gp = Omega^infty(K)}
  Let $\sC_{/S} \sub \Aff_{/S}$ be any admissible subcategory.
  Denote by $X_S$ the presheaf on $\sC_{/S}$ given by the assignment $\Spec(R) \mapsto X(R)$, and by $(X_S)^\gp$ its group completion.
  Then there is a canonical morphism of $\sC$-fibred $\Einfty$-groups
  \begin{equation*}
    (X_S)^\gp \to \Omega^\infty(\K|_{\sC})
  \end{equation*}
  which induces an isomorphism $\LZar (X_S)^\gp \simeq \Omega^\infty(\K|_{\sC})$.
\end{prop}

\begin{proof}
  Let $X'_S$ denote the presheaf $\Spec(R) \mapsto (\Mod_R^\proj)^\simeq$.
  Then by \cite[Lect.~19, Thm.~5]{LurieKLectures} there is a canonical isomorphism $(X'_S)^\gp \simeq \Omega^\infty(\K|_{\sC})$.
  Therefore it will suffice to show that the monomorphism $X_S \hook X'_S$ induces an effective epimorphism of Zariski sheaves $\LZar(X_S) \to X'_S$ (see \cite[Ex.~5.2.8.16]{HTT}).
  This is clear since every finitely generated projective $R$-module is Zariski-locally free.
\end{proof}

\begin{rem}\label{lem:Li_! connective}
  Note that the functors $\LNis\iota_!$ and $\L\iota_!$ preserve connective objects.
  Indeed, the essential images of the fully faithful functors
  \begin{align*}
    \MotSpc(\sA_{/S})_{\Spt_{\geq 0}} \hook \MotSpc(\sA_{/S})_{\Spt},\\
    \MotSpc(\sB_{/S})_{\Spt_{\geq 0}} \hook \MotSpc(\sB_{/S})_{\Spt},
  \end{align*}
  are generated under colimits by objects of the form $\Sigma^{\infty}_+(X)$ with $X \in \sA_{/S}$, resp. $X \in \sB_{/S}$.
\end{rem}

\begin{proof}[Proof of \propref{prop:unstable comparison}]
  Since both source and target are connective, it will suffice to show the claim for the underlying $\sA$-fibred $\Einfty$-group $\Omega^\infty(\K|_{\sA})$.
  Since each of the functors $\LNis$, $\iota_!$ and $\iota^*$ commutes with colimits and finite products, and hence with group completion (see e.g. \cite[Lem.~5.5]{HoyoisCdh}), we reduce using \propref{prop:X^gp = Omega^infty(K)} to showing the analogous claim for the $\sA$-fibred $\Einfty$-group $X_S$, i.e., that the canonical morphism
  \[
  \L w_! (X_S) \simeq X^\cl_{S_\cl}
  \]
  is invertible, where the right-hand side is the construction analogous to \constrref{constr:X} in classical algebraic geometry.
  The claim now follows from the fact that $w : \sA_{/S} \to \sA^\cl_{/S_\cl}$ sends the spectral affine schemes\footnote{%
    Since $\GL_{n,S}$ are Zariski-open inside spectral affine spaces, they belong not only to $\Sm_{/S}$ but even to $\sA^0_{/S}$ (\examref{exam:minimal adm}) and hence to any narrow $\sA_{/S}$.
  } $\GL_{n,S}$ to their classical counterparts.
\end{proof}

\begin{proof}[Proof of \propref{prop:unstable K Sm-fibred}]\leavevmode
  As above, we reduce to the analogous claim for the $\sB$-fibred $\Einfty$-group $X_S$.
  This follows from the fact that the classifying spaces $\BGL_{n,S}$ are colimits of finite products of the smooth spectral schemes $\GL_{n,S}$.
\end{proof}


\ssec{Bott periodicity}
\label{ssec:kh/bott}

We keep the notation of \ref{notat:KHlafds}.
In this subsection we will show that $\KH|_{\sB}$ can be described as the Bott periodization of the motivic localization of $\K|_{\sB}$, and similarly for $\KH|_{\sA}$.
This will allow us to prove a nonconnective analogue of \propref{prop:unstable K Sm-fibred}, see \corref{cor:K Sm-fibred}.

\begin{constr}\label{constr:Bott}
  Let $R$ be a \cnEring.
  Denote by $R\{T\}$ the free $\Einfty$-$R$-algebra on one generator $T$ and by $R\{T,T^{-1}\}$ the localization away from $T \in \pi_0(R\{T\}) \simeq \pi_0(R)[T]$.
  Just as in \remref{rem:Bott class}, the automorphism of $R\{T,T^{-1}\}$ given by multiplication by $T$ induces a canonical element $b \in \K_1(R\{T,T^{-1}\})$ which we also call the Bott class.
  Again by \cite[Lect.~20, Cor.~4]{LurieKLectures} there is a canonical isomorphism $\K_1(R\{T,T^{-1}\}) \simeq \K_1(\pi_0(R)[T,T^{-1}])$ under which $b$ corresponds to the usual Bott class.
  In particular, the canonical bijection $\K_1(R\{T,T^{-1}\}) \to \K_1(R[T,T^{-1}])$ (induced by $\varepsilon$, see \remref{rem:A^n_flat to A^n}) sends $b$ to $b$.
\end{constr}

\begin{rem}
  The Bott class may be regarded as a morphism of $\sB$-fibred $S^1$-spectra
  \[
    b : \Sigma^\infty(\G_{m,S},1)[1] \to \K|_{\sB}.
  \]
\end{rem}

To simplify notation, we set $\T_S := \Sigma^\infty(\G_{m,S},1)[1]$.

\begin{defn}\label{constr:Bott periodic}
  Recall that the $\sB$-fibred $S^1$-spectrum $\K|_{\sB}$ admits an $\Einfty$-ring structure, induced by the symmetric monoidal structure on perfect complexes.
  We say that a $\K|_{\sB}$-module $\sF$ is \emph{Bott-periodic} if the canonical morphism
  \[
    b : \sF \to \uHom(\T_S, \sF)
  \]
  induced by the Bott class (via the action of $\K|_{\sB}$ on $\sF$) is invertible.
  We define Bott-periodic $\K|_{\sA}$-modules similarly.
  Note that the full subcategory spanned by Bott-periodic $\K|_{\sB}$-modules (resp. $\K|_{\sA}$-modules) is a left localization; we refer to the left adjoint $\Q$ as \emph{Bott periodization}.
\end{defn}

\begin{thm}\label{thm:KH as Bott periodization}
The canonical morphisms
  \begin{align*}
    \K|_{\sB} \to \KH|_{\sB}\\
    \K|_{\sA} \to \KH|_{\sA}
  \end{align*}
induce isomorphisms
\begin{align*}
  \Q\big(\L (\K|_{\sB})\big) &\simeq \KH|_{\sB},\\
  \Q\big(\L (\K|_{\sA})\big) &\simeq \KH|_{\sA}\\
\end{align*}
of Bott-periodic motivic fibred $S^1$-spectra.
\end{thm}

The following description of Bott periodization will be used in the proof of \thmref{thm:KH as Bott periodization}.

\begin{rem}\label{rem:Q_b}
  Explicitly, the Bott periodization of a \emph{motivic} $\K|_{\sB}$-module $\sF$ can be computed as the colimit of the tower
  \[
    \sF
      \xrightarrow{b} \uHom(\T_S, \sF)
      \xrightarrow{b} \uHom(\T_S^{\otimes 2}, \sF)
      \xrightarrow{b} \cdots
  \]
  according to Theorem~3.8 and the proof of Lemma~4.9 of \cite{HoyoisCdh}.
  Moreover, if $\sF$ is nil-local, then we may replace $b$ by the Bott class in $\G_{m,S}^\flat$ (\remref{rem:Bott class}).
\end{rem}

\begin{rem}\label{rem:i^* Bott}
  Since the functors $\L\iota_!$ and $\iota^*$ are symmetric monoidal (Remarks~\ref{rem:i^*Hom} and \ref{rem:tensor V}), they extend to an adjunction
  \[
    \L\iota_! : \Mod_{\L(\K|_{\sA})}(\MotSpc(\sA_{/S})_\Spt)
      \longrightleftarrows \Mod_{\L(\K|_{\sB})}(\MotSpc(\sB_{/S})_\Spt) : \iota^*.
  \]
  Since $\T_S$ belongs to the essential image of $\L\iota_!$ (as $\G_{m,S}$ belongs to $\sA_{/S}$), it follows from \remref{rem:i^*Hom}(iii) and the fact that $\iota^*$ preserves colimits that $\iota^*$ preserves Bott-periodic objects and commutes with the Bott periodization functor $\Q$.
  Its left adjoint on Bott-periodic objects is given by $\sF \mapsto \Q(\L\iota_!(\sF))$ and preserves $\Q$-equivalences.
\end{rem}

\begin{proof}[Proof of \thmref{thm:KH as Bott periodization}]
  Since $\iota^*$ commutes with $\Q$ (\remref{rem:Q_b}), it will suffice to consider the first map.
  Just as in the classical setting \cite[Prop.~2.10]{cisinski2013descente}, one observes that up to $\A^{1,\flat}$-localization, hence also up to $\A^1$-localization by \lemref{lem:A^1 vs A^1_flat-invariance}, the Bass construction (\constrref{constr:Bass}) simplifies to give the formula
  \[
    \KH|_{\sB} \simeq \colim \left(
      \L(\K|_{\sB})
      \xrightarrow{b} \uHom\left(\Sigma^\infty(\G^\flat_{m,S},1)[1], \L(\K|_{\sB})\right)
      \xrightarrow{b} \cdots
    \right),
  \]
  where the maps are induced by the Bott class $b\in\K_1(\bG^\flat_{m,S})$ (\remref{rem:Bott class}).
  But by \remref{rem:Q_b} and \corref{cor:LK|B nil-local}, this is isomorphic to the Bott periodization $\Q(\L(\K|_{\sB}))$ with respect to $b\in\K_1(\bG_{m,S})$.
\end{proof}

Using \thmref{thm:KH as Bott periodization} we may deduce the following $S^1$-stable analogue of \propref{prop:unstable K Sm-fibred}, which holds up to $\A^1$-homotopy and Bott periodization:

\begin{cor}\label{cor:K Sm-fibred}
  The canonical morphism
  \begin{equation*}
    \L \iota_! (\KH|_{\sA}) \simeq \L \iota_!\iota^* (\KH|_{\sB}) \to \KH|_{\sB}
  \end{equation*}
  induces an isomorphism
  \[
    \Q\big(\L\iota_! (\KH|_{\sA})\big) \simeq \KH|_{\sB}
  \]
  of Bott-periodic motivic $\sB$-fibred $S^1$-spectra.
\end{cor}
\begin{proof}
  Follows immediately from \thmref{thm:KH as Bott periodization} and \remref{rem:i^* Bott}.
\end{proof}

\begin{rem}
In the statement of \corref{cor:K Sm-fibred}, we can replace the source with $\L \iota_! (\KB|_{\sA})$.
That is, the canonical map
  \begin{equation*}
    \Q(\L\iota_! (\KB|_{\sA})) \to \KH|_{\sB}
  \end{equation*}
is also invertible.
This follows from \propref{prop:KH motivic} and the fact that $\iota_!$ preserves motivic equivalences (\lemref{lem:iota and local equivs}).
\end{rem}


\ssec{Proof of \thmref{thm:v_!KH=KH}}
\label{ssec:kh/coda}

The only remaining ingredient is the behaviour of Bott periodization under the equivalence
\[
  \L w_! : \MotSpc(\sA_{/S}) \to \MotSpc(\sA^\cl_{/S_\cl})
\]
of \thmref{thm:intro/comp1}.
But the fact that it commutes with internal homs (\remref{rem:functoriality of nil-localization}\itemref{item:aufpasdmf}) immediately implies that we have
\[
  \L w_! (\Q(\sF)) \simeq \Q^\cl(\L w_!(\sF))
\]
in $\MotSpc(\sA^\cl_{/S_\cl})$ for every $\L(\K)$-module $\sF$ in $\MotSpc(\sA_{/S})$.
Here $\Q^\cl$ denotes Bott periodization of a $\SmCl$-fibred motivic space with the classical Bott element.

Now consider the commutative diagram in $\MotSpc(\sA^\cl_{/S_\cl})$
\[ \begin{tikzcd}
  \L w_! \left(\Q(\L(\K|_{\sA}))\right) \ar{r}\ar{d}
    & \L w_!(\KH|_{\sA}) \ar{d}
  \\
  \Q(\L(\K|_{\sA^\cl})) \ar{r}
    & \KH^\cl|_{\sA^\cl}.
\end{tikzcd} \]
The assertion of \thmref{thm:v_!KH=KH} is that the right-hand vertical arrow is invertible.
Since the horizontal arrows are isomorphisms by \thmref{thm:KH as Bott periodization} (and its classical analogue), it will suffice to demonstrate the invertibility of the left-hand vertical arrow.
Since $\L w_!$ commutes with $\Q$, this is identified with the Bott periodization of the canonical morphism
\[ \L w_!(\K|_{\sA}) \simeq \L w_!(\L (\K|_{\sA})) \to \L(\K|_{\sA^\cl}), \]
which is invertible by \propref{prop:unstable comparison}.

\appendix

\bibliographystyle{alphamod}

\renc{\baselinestretch}{1.0}
{\small
\bibliography{references}
}

\end{document}